\documentclass[twoside,12pt]{amsart}
\usepackage{amsmath}
\usepackage{amsthm}
\usepackage{latexsym}
\usepackage{amssymb}

\def\bdm{\begin{displaymath}}
\def\edm{\end{displaymath}}

\def\bba{{\mathbf a}}
\def\bbb{{\mathbf b}}
\def\bbA{{\mathbf A}}
\def\bbB{{\mathbf B}}
\def\Aeta{{\mathbf A}^{\eta}}
\def\Aetaj{{A}_j^{\eta}}
\def\Beta{{\mathbf B}^{\eta}}
\def\Betajk{{B}_{jk}^{\eta}}

\newtheorem{thm}{Theorem}[section]
\newtheorem{avg}{Averaging Lemma}[section]

\newtheorem{cor}[thm]{Corollary}

\theoremstyle{definition}

\theoremstyle{remark}
\newtheorem{rem}[thm]{Remark}
\newtheorem{remark}[thm]{Remark}

\numberwithin{equation}{section}

\def\cal{\mathcal}

\def\ixi{i\xi}
\def\xip{\xi'}
\def\ixip{\ixi'}
\def\tip{\tau'}

\def\cutL1{\cap L^1}

\def\supa{\mathop{\sup}}
\def\supxi{\supa_{|\xi|=1}}    

\newcommand{\ff}{{\cal F}}

\newcommand{\dd}{{{\cal D}}}

\def\bi{\begin{itemize}}
\def\bs{\begin{split}}
\def\es{\end{split}}
\def\ba{\begin{align}}
\def\bas{\begin{align*}}
\def\ea{\end{align}}
\def\eas{\end{align*}}

\def\trace{{\operatorname{trace}}}

\def\sgn{{\operatorname{sgn}}}

\def\J{{\cal J}}
\def\eps{\varepsilon}

\def\lesssim{\stackrel{{}_<}{{}_\sim}}
\def\gtrsim{\stackrel{{}_>}{{}_\sim}}

\def\emph#1{{\it #1}}
\def\textbf#1{{\bf #1}}

\def\myM{{\cal M}}
\def\mym{N}		
\def\vx{x,v}
\def\tx{t,x}
\def\txv{\tx,v}
\def\vnab{\nabla_x,v}
\def\vxi{\xi,v}
\def\vixi{\ixi,v}

\def\vixip{\ixip,v}
\def\op{{\mathcal L}}		
\newcommand{\symb}{\op(\ixip,v)}
\def\myo{\omega_{\op}(\dyM;\delta)}
\def\myO{\Omega_{\op}(\xi;\delta)}
\def\dyM{J}				

\newcommand{\R}{\mathbb{R}}

\newcommand{\C}{\mathbb{C}}
\newcommand{\Z}{\mathbb{Z}}
\newcommand{\loc}{\mbox{\scriptsize loc}}

\newcommand{\Rd}{\R^d}
\newcommand{\Rt}{(0,\infty)}
\newcommand{\Rdx}{\R^d_x}
\newcommand{\Rtdx}{\R_t\times\Rdx}

\newcommand{\Rv}{\R_v}
\newcommand{\Rvdx}{\Rdx\times\Rv}
\newcommand{\Rtdxv}{\Rtdx\times\Rv}
\newcommand{\Wsp}[1]{W^{\sigma,p}_{loc}({#1})}
\newcommand{\Wsr}[1]{W^{s,r}_{loc}({#1})}
\newcommand{\Lqloc}{L^q_{loc}}
\newcommand{\Lploc}{L^p_{loc}}
\newcommand{\stheta}{(1-\theta)\sigma+\theta(k-\eta)}

  \newtheorem{lemma}{Lemma}[section]

  \newtheorem{definition}{Definition}[section]

\begin{document}

\title[Velocity averaging and regularity of quasilinear PDEs]{Velocity averaging, kinetic formulations\\ and regularizing effects in quasilinear PDEs}



\author{Eitan Tadmor}
\address{Department of Mathematics, Center of Scientific Computation And Mathematical Modeling (CSCAMM)
and Institute for Physical Science and Technology (IPST), University of Maryland, MD 20742.}
\curraddr{} \email{tadmor@cscamm.umd.edu}

\author{Terence Tao}
\address{Department of Mathematics, UCLA, Los Angeles, CA 90095}
\curraddr{} \email{tao@math.ucla.edu}

\thanks{{\tiny Part of this research was carried out while E.T. was visiting the Weizmann Institute for Science and it is a pleasure to thank the Faculty of Mathematics and Computer Science for their hospitality. Research was supported in part  by NSF DMS \#DMS04-07704 and ONR \#N00014-91-J-1076 (E.T.)
and by the  David and Lucille Packard foundation (T.T). }}
\date{\today}



\begin{abstract}
We prove new velocity averaging results for second-order multidimensional equations of the general form, 
$\op(\nabla_x,v)f(x,v)=g(x,v)$ where  $\op(\nabla_x,v):=\bba(v)\cdot\nabla_x-\nabla_x^\top\cdot\bbb(v)\nabla_x$. These results quantify the Sobolev regularity of the averages, $\int_vf(x,v)\phi(v)dv$, in terms of the non-degeneracy of the set $\{v\!: |\op(\ixi,v)|\leq \delta\}$ and the mere integrability of the data, $(f,g)\in (L^p_{x,v},L^q_{x,v})$.  Velocity averaging is then used to study the \emph{regularizing effect} in quasilinear second-order equations, $\op(\nabla_x,\rho)\rho=S(\rho)$ using their underlying kinetic formulations, $\op(\nabla_x,v)\chi_\rho=g_{{}_S}$. In particular, we improve previous regularity statements for nonlinear conservation laws, and we derive completely new regularity results for convection-diffusion and elliptic equations driven by degenerate, non-isotropic diffusion.       
\end{abstract}

\maketitle

\bigskip
\centerline{To Peter Lax and Louis Nirenberg on their $80^{th}$ birthday}
\centerline{With friendship and admiration}

{\small \tableofcontents}

\section{Introduction}
We study the regularity of solutions to multidimensional quasilinear scalar equations of the form
\begin{equation}\label{eq:class}
\sum_{j=0}^d\frac{\partial}{\partial x_j}A_j(\rho)
-\sum_{j,k=1}^d \frac{\partial^2}{\partial x_j\partial x_k}B_{jk}(\rho) = S(\rho) 
\end{equation}
where $\rho: \R^{1+d} \to \R$ is the unknown field and $A_j$, $B_{jk}$, $S$ are given functions from $\R$ to $\R$.

This class of equations governs time-dependent solutions, $\rho(t\!=\!x_0,x)$, of nonlinear conservation laws  where $A_0(\rho)=\rho$ and $B_{jk}\equiv 0$, time-dependent solutions of degenerate diffusion  and convection-diffusion equations where $\{B'_{jk}\}\geq 0$, and spatial solutions, $\rho(x)$, of degenerate elliptic equations  where $A_j\equiv0$.
The notion of solution should be interpreted here in an appropriate weak sense, since we focus our attention on 
the degenerate diffusion case, which is too weak to enforce the smoothness required for a notion of a strong solution.  Instead, a common feature of such problems is the (limited) regularity of their solutions, which is dictated by the \emph{nonlinearity} of the governing equations. A prototype example is provided by discontinuous solutions of nonlinear conservation laws. In  \cite{LPT94a}, Lions, Perthame \& Tadmor have shown  that entropy solutions of such laws admit a regularizing effect of a fractional order,  dictated  by the order of non-degeneracy of the equations. In this paper we extend this result for the general class of second-order equations (\ref{eq:class}). In particular, we improve the \cite{LPT94a}-regularity statement for nonlinear conservation laws, and we derive completely new regularity results for convection-diffusion and elliptic equations driven by degenerate, non-isotropic diffusion.       

The derivation of these regularity results employs a \emph{kinetic formulation} of (\ref{eq:class}). To describe this formulation let us proceed formally,  seeking an equation which governs the indicator function, $\chi_{\rho(x)}(v) := \sgn(v)(|\rho|-|v|)_+$ associated with $\rho$, and which depends on an auxiliary  velocity variable $v\in\R$, borrowing the terminology from the classical kinetic framework.   To this end, we consider the distribution $g=g(x,v)$, defined via its velocity derivative $\partial_v g$ using the formula

\begin{equation}\label{eq:kin-class}
\partial_v g(x,v):=\Big(\bba(v)\cdot\nabla_x -\nabla_x^\top\cdot\bbb(v)\nabla_x + S(v)\partial_v \Big)\chi_{\rho(x)}(v), \quad \bba_j:=A_j', \ \bbb_{jk}:=B'_{jk}\geq0.
\end{equation}
Observe that the nonlinear quantities $\Phi(\rho)$ can be expressed as the $v$-moments of $\chi_{\rho}$, 
$\Phi(\rho) \equiv \int_v \Phi'(v)\chi_{\rho}(v)dv, \ \Phi(0)=0$. Therefore, by velocity averaging
of (\ref{eq:kin-class}) we recover (\ref{eq:class}). Moreover, for a proper notion of weak solution $\rho$, one augments (\ref{eq:class}) with additional conditions on the behavior of $\Phi(\rho)$ for a large enough family of entropies $\Phi$'s. These additional entropy conditions imply that $g$ is in fact a positive distribution, $g=m \in {\cal M}^+$, measuring the entropy dissipation of the nonlinear equation.
We arrive at the kinetic formulation of (\ref{eq:class})

\begin{equation}\label{eq:for-class}
\op(\nabla_x,v)f(x,v)=\partial_v m-S(v)\partial_v\chi_{\rho(x)}(v), \quad  f=\chi_\rho, \ m\in {\cal M}^+
\end{equation}  
where $\op$ is identified with the linear symbol 
$\op(\ixi,v):=\bba(v)\cdot\ixi+\langle\bbb(v)\xi,\xi\rangle$.
We recall that $\rho(x)$ itself can be recovered by velocity averaging of $f(x,v)=\chi_\rho(v)$, via the identity
$\rho(x) = \int f(x,v)\ dv$.
In Section \ref{sec:avg} we discuss the regularity gained by such velocity averaging.
There is a relatively short yet intense history of such regularity results, commonly known 
as `velocity averaging lemmata'. We mention the early works of \cite{GLPS88}, \cite{DLM91} and  their applications, 
in the context of nonlinear  conservation laws, in \cite{LPT94a, LPT94b,LPS96, JP02}; a detailed list of references can be found in \cite{Pe02} and is revisited in Section \ref{sec:avg} below. Almost all previous averaging results  dealing with (\ref{eq:for-class}) were restricted to    
first-order transport equation, $\deg\op(\ixi,\cdot)=1$.
In Section \ref{subsec:hom} we present an extension to general symbols $\op$'s which satisfy the so called \emph{truncation property}; luckily --- as shown in Section \ref{subsec:fssymbols} below, all $\op$'s with $\deg\op(\ixi,\cdot)\leq 2$ satisfy this property. If  the problem is \emph{non-degenerate} in the  sense that\footnote{We use $X \lesssim Y$ to denote the estimate $X \leq CY$ where $C$ is a constant which can depend on exponents such as $\alpha, \beta, p$ and on symbols such as $\op$, $A_j$, $B_{jk}$ but is independent of fields such as $\rho$, co-ordinates such as $x,t,v$, and scale parameters such as $\delta$.  We use $X \sim Y$ to denote the assertion that $X \lesssim Y \lesssim X$.}
\[
\exists \alpha \in (0,1), \ \beta>0 \ \ \hbox{s.t.} \ \ \sup_{|\xi|\sim J} \big|\big\{v: \ |\op(\vixi)|\leq \delta\big\}\big| \lesssim \Big(\frac{\delta}{J^\beta}\Big)^\alpha,
\]
then  by velocity averaging of the kinetic solution $f$  in  (\ref{eq:for-class}), $\rho(x)=\int_v f(x,v)dv$ has a  $W^{s,r}$-regularity of order $s<\beta\alpha/(3\alpha+2)$ with an appropriate $r=r_s>1$; consult (\ref{eq:opts}) below. In the particular case of non-degenerate homogeneous symbols of order $k$, 
\[
\exists \alpha \in (0,1)\ \ \hbox{s.t.} \ \ \sup_{|\xi|=1} \big|\big\{v: \ |\op(\vixi)|\leq \delta\big\}\big| \lesssim \delta^\alpha,
\]
we may take $\beta=k$ and  velocity averaging implies   a $W^{s,r}_{loc}$-regularity exponent $s<k\alpha/(3\alpha+2)$. The main results are summarized in averaging lemmata \ref{thm:improved-averaging} and \ref{avg}.

In Section \ref{sec:hyp} we turn to the first application of these averaging results in the context of nonlinear conservation laws, $\rho_t+\nabla\cdot A(\rho)=0$, subject to $L^\infty$-initial data $\rho_0$. If the equation is non-degenerate of order $\alpha$ in the sense that
\[
\big|\left\{v: |\tau+\bba(v)\cdot\xi|\leq\delta\right\}\big| \lesssim \delta^\alpha \ \  \ \hbox{and} \ \ \sup_{\left\{v: \ |\tau+\bba(v)\cdot\xi|\leq\delta\right\}}|\bba'(v)\cdot\xi|\lesssim \delta^{1-\alpha}, \quad \forall \ \tau^2+|\xi|^2=1,
\]
then for $t>0$, the entropy solution $\rho(t,x)$ gains Sobolev regularity  $\rho(t,\cdot)\in W^{s,1}_{\loc}(\Rdx)$ of order $s<s_\alpha=\alpha/(2\alpha+1)$. This improves the Sobolev-regularity exponent of order $\alpha/(\alpha+2)$ derived at \cite{LPT94a} (while facing the same barrier of $s_1=1/3$ discussed in \cite{DLW05}).

Section \ref{sec:par} is devoted to convection-diffusion equations. We begin, in Section \ref{subsec:deg}, with second-order degenerate diffusion $\rho_t -\sum \partial^2_{x_jx_k}B_{jk}(\rho)=0$. The emphasis here is on non-isotropic diffusion, beyond the prototype case of the porous medium equation (which corresponds to the case when $B_{jk}$ is a scalar multiple of the identity,
$B_{jk} = B \delta_{jk}$). The regularizing effect is determined by the smallest non-zero eigenvalue $\lambda(v)=\lambda(\bbb(v))\equiv \!\!\!\!\!\!/ \   0$ of $\bbb(v):=B'(v)$, so that
\begin{equation}\label{eq:indeg}
\big|\left\{v: 0\leq\lambda(v)\leq \delta\right\}\big| \lesssim \delta^\alpha \ \  \ \hbox{and} \ \ \sup_{\left\{v:\ \lambda(v)\leq\delta\right\}}|\langle\bbb'(v)\xi,\xi\rangle|\lesssim \delta^{1-\alpha}, \quad \forall \ |\xi|=1.
\end{equation}
Staring with initial conditions $\rho_0\in L^\infty$, then the corresponding kinetic solution $\rho(t>0,\cdot)$ gains  $W^{s,1}_{\loc}$-regularity of order $s<2\alpha/(2\alpha+1)$. 
In Section \ref{subsec:cd} we take into account the  additional effect of  nonlinear convection. The resulting convection-diffusion equations, coupling  degenerate and possibly non-isotropic diffusion with non-convex convection governing capillarity effects, are found in a variety of applications.
Consider the prototype one-dimensional case  
\[
\rho_t+ A(\rho)_x- B(\rho)_{xx}=0,  \quad A'(\rho)\sim \rho^\ell, B'(\rho)\sim |\rho|^n.
\]
The regularizing effect is dictated by the strength of the degenerate  diffusion vs. the convective degeneracy. If $n\leq \ell$ we find
$W^{s,1}_{loc}$-regularity of order $s<2/(n+2)$, which is the same Sobolev-regularity exponent we find with the `purely diffusive' porous-medium equation, i.e., when $A=0$. On the other hand, if the diffusion is too weak so that $n\geq 2\ell$, we then conclude with a Sobolev-regularity exponent of order $s<1/(\ell+2)$, which is dominated by the convective part of the equation. In Section \ref{subsec:cd} we present similar results for
multidimensional convection-diffusion equations with increasing degree of degeneracy. In particular, consider the two-dimensional equation
\[
\rho_t+(\partial_{x_1}+\partial_{x_2})A(\rho)-(\partial_{x_1}-\partial_{x_2})^2B(\rho)=0.
\]
If we set $A=0$, the equation has a strong, rank-one parabolic degeneracy with no regularizing effect coming from its purely diffusion part,
since $\langle\bbb(v)\xi,\xi\rangle\equiv 0, \ \forall \xi_1-\xi_2=0$, indicating the persistence of steady oscillations,
$\rho_0(x+y)$; moreover, if $B=0$ then the equation has no regularization coming from its purely convection part, since $\bba(v)\cdot \xi\equiv 0, \ \forall \xi_1+\xi_2=0$ indicates the persistence of steady oscillations $\rho_0(x-y)$. 
Nevertheless, the combined \emph{convection-diffusion} with $A(\rho)\sim \rho^{\ell+1}$ and $B(\rho)\sim |\rho|^n\rho$  does have  $W^{s,1}$-regularizing effect of order $s<6/(2+2n-\ell)$ for $n\geq 2\ell$, consult corollary \ref{cor:mix} below.

Finally, in Section \ref{sec:ellip} we consider degenerate elliptic equations, 
$-\sum \partial^2_{x_jx_k}B_{jk}(\rho)=S(\rho)$. Assuming that the non-degeneracy condition (\ref{eq:indeg}) holds, then the kinetic formulation of bounded solutions for such equations, $\rho\in W^{s,1}(D)$, have interior regularity of order $s<\min(\alpha, 2\alpha/(2\alpha+1))$.
We conclude by noting that it is possible to adapt our arguments in more general setups, for example, when a suitable
source term $S(\rho)$ it added to the time dependent problem, when dealing with degenerate \emph{temporal} fluxes, $\partial_t A_0(\rho), \ A'_0\geq0$,  or when lower-order convective terms, $\nabla_x\cdot A(\rho)$, are added in the elliptic case.  

\section{The averaging lemma}\label{sec:avg}

We are concerned with the regularity of averages of solutions for differential equations of the form
\begin{equation}\label{eq:transport}
 \op(\vnab)f=\Lambda_x^\eta\partial^\mym_v g, \quad  \Lambda_x:=(-\Delta_x^2)^{1/2}.
\end{equation}
Here, $f=f(\vx)\in\Wsp{\Rvdx}$ and $g=g(\vx)\in\Lqloc(\Rvdx)$ are real-valued functions of the spatial variables 
$x=(x_1,\ldots,x_d)\in \Rd$ and an additional parameter $v\in\R$, called \emph{velocity} by analogy with the kinetic framework, and $\op(\vnab)$ is a differential operator on $\Rdx$ of order $\leq k$, whose  coefficients are smooth functions of $v$. 
 
The velocity averaging lemma asserts that if  $\op(\cdot,v)$ is \emph{nondegenerate} in the sense  that its null set is sufficiently small ---  to be made precise below,
 then the $v$-moments of $f(x,\cdot)$,
\bdm
\overline{f}(x):= \int_v f(\vx)\phi(v)dv, \qquad \phi \in C_0^\infty,
\edm
are smoother than the usual regularity associated with the 
data of $f(x,\cdot)$ and $g(x,\cdot)$. That is, by averaging over the so-called microscopic $v$-variable, 
there is a  \emph{gain} of regularity in the macroscopic $x$-variables.
There is a relatively short yet intense history of such results,  
motivated by kinetic models such as Boltzmann, Vlasov, radiative transfer and similar
equations where the $v$-moments of $f$ represent macroscopic quantities of interest. 
We refer to the early works of Agoshkov \cite{Ag84} and Golse, Lions, Perthame and Sentis, \cite{GPS85, GLPS88}
 treating first-order transport operators with $f,g$ integrability of order $p=q>1$. The  work of DiPerna, Lions and Meyer, \cite{DLM91}, provided the first treatment of the general  case $p\neq q$, followed by B\'{e}zard, \cite{Be94},
their optimality in \cite{Li95} and an optimal Besov regularity result of DeVore and Petrova, \cite{DVP00}.

Extensions to more general streaming operators were treated by DiPerna and Lions in \cite{DL89a, DL89b} with applications to Boltzmann and Vlasov-Maxwell equations, and by Lions, Perthame and Tadmor in \cite{LPT94a} with applications
to nonlinear conservation laws and related parabolic equations.  G\'{e}rard, \cite{Ge90}  together with Golse \cite{GG92}   provided an $L^2$-treatment of general differential operators.
A different line of extensions consists of  velocity  lemmata which take into account different orders of integrability
in $x$ and $v$,  leading to sharper  velocity  regularization results sought in various application. We refer to  the results in \cite{DLM91} for $f \in L^{p_1}(\Rv, L^{p_2}(\Rdx)), g \in L^{q_1}(\Rv, L^{q_2}(\Rdx))$ and to  Jabin and Vega \cite{JV03} for  $f \in W^{N_1,p_1}(\Rv, L^{p_2}(\Rdx)), g \in W^{N_2,q_1}(\Rv, L^{q_2}(\Rdx))$. Westdickenberg, \cite{We02}, analyzed a general case of the form  $f \in B^{\sigma}_{p_1,p_2}(\Rdx, L^{r_1}(\Rv)), g \in B^{\sigma}_{q_1,q_2}(\Rdx, L^{r_2}(\Rv))$.
Jabin, Perthame and Vega, \cite{JP02}, \cite{JV04} used a mixed integrability of $f \in L^p(\Rdx, W^{N_1,p}(\Rv)), g \in L^q(\Rdx, W^{N_2,q}(\Rv))$ to improve the regularizing results for nonlinear conservation laws \cite{LPT94a}, \cite{LPT94b}, Ginzburg-Landau, and other nonlinear models. Golse and Saint-Raymond, \cite{GSR02}  have shown that a minimal requirement of \emph{equi-integrability}, say   $f\in L^1(\Rdx,L^\Phi(\Rv)), g \in L^1(\Rdx \times \Rv)$ measured in Orlicz space $L^\Phi$ with super-linear $\Phi$,  is sufficient for relative compactness of the averages $\overline{f}$, which otherwise might fail for mere $L^1$-integrability, \cite{GLPS88}. 

The derivation of velocity averaging in the above works was accomplished by  various methods. The main approach, which we use below, is based on decomposition in Fourier space, carefully tracking 
$\widehat{f}(\xi,v)$ in the ``elliptic" region where $\{v  : \ \op(\vixi) \neq 0 \}$, and 
the complement region which is made sufficiently small by a non-degeneracy assumption. Other approaches include the use  microlocal defect measures and H-measures  \cite{Ge90, Ta90}, wavelet decomposition \cite{DVP00} and ``real-space methods" --- in time \cite{Va99, BD99}, and in space-time using Radon transform \cite{Ch00, We02}, X-transform  \cite{JV03,JV04} and duality-based dispersion estimates \cite{GSR02}. Almost all of these results  are devoted to the phenomena of velocity averaging in the context of transport equations, $k=1$,  

Our study of velocity  averaging applies to a large class of  $\op$'s, satisfying the so-called \emph{truncation property}:
in Section \ref{subsec:fssymbols} we show that all  $\op$'s of order $k \leq 2$ satisfy this truncation property. In particular, we improve the regularity statement for first-order velocity averaging and extend the various velocity averaging results of the works above from first-order transport to general second-order transport-diffusion and elliptic equations. The results are summarized in the averaging lemmata \ref{thm:homogeneous-averaging} and \ref{thm:improved-averaging} for homogeneous symbols and in averaging lemma \ref{avg}  for general, truncation-property-satisfying $\op$'s. Our derivation is carried out in Fourier space using Littlewood-Paley decompositions of $f\in\Wsp{\Rvdx}$ and $g\in\Lqloc(\Rvdx)$. To avoid an overload of indices, we leave for future work possible extensions for more general   data with mixed $(x,v)$-integrability of $f$ and $g$.

\subsection{The truncation property}\label{subsec:trunc}

We now come to a fundamental definition.

\begin{definition}
 Let $m(\xi)$ be a complex-valued Fourier multiplier.  We say that $m$ has the \emph{truncation property} if, for any locally-supported bump function $\psi$ on $\C$ and any $1 <  p < \infty$, the multiplier with symbol $\psi( m(\xi)/\delta )$ is an  $L^p$-multiplier  uniformly in $\delta > 0$, that is, its $L^p$-multiplier norm depends solely on the support and $C^\ell$ size of $\psi$ (for some large $\ell$ which may depend on $m$) but otherwise is independent of $\delta$.
\end{definition}

\smallskip\noindent
In Section \ref{subsec:fssymbols} we will describe some examples of multipliers with this truncation property.
Equipped with this notion of a truncation property, we turn to discuss the $L^p$-size of \emph{parameterized multipliers}.  Let
\begin{align*}
\myM_\psi f(\vx)&:= {\ff}_x^{-1}\psi\Big(\frac{m(\vxi)}{\delta}\Big)\widehat{f}(\vxi), \hbox{ where }
\widehat{f}(\vxi)={\ff}_x f(\vxi) := \int_{\R^d} e^{-2\pi i \xi \cdot x} f(x,v)\ dx
\end{align*}
denote the operator associated with the complex-valued multiplier $m(\cdot,v)$, truncated at level $\delta$. 
Our derivation of the averaging lemmata below  is based on the following straightforward
estimate for such parameterized multipliers.

\begin{lemma}[Basic $L^p$ estimate]\label{basic-lp}
Let $I$ be a finite interval, $I \subset \Rv$ and assume $m(\xi,v)$ satisfies the truncation property uniformly in $v\in I$. Let $1<p \leq 2$.  Let $\overline{\myM_\psi}$ denote the velocity-averaged Fourier multiplier
$$ \overline{\myM_\psi} f(x) := \int_I M_\psi f(x,v)\ dv
= \int_I \ff_x^{-1}\psi(m(\xi,v)/\delta)\ff_x f(x,v)\ dv.$$ 
For each $\xi \in\R^d$ and $\delta > 0$ let $\Omega_m(\xi;\delta) \subset I$ be the velocity set
$$ \Omega_m(\xi;\delta):=\Big\{v \in I \ : \ \frac{m(\xi,v)}{\delta} \in supp \ \psi \Big\}.$$
Then we have the $L^p$ multiplier estimate
\begin{equation} \label{eq:averaged-multiplier}
\|\overline{\myM_\psi f}(x)\|_{L^p(\Rdx)}  \lesssim   \sup_\xi|\Omega_m(\xi;\delta)|^{\frac{1}{p'}}\cdot\|f\|_{L^p(\Rvdx)}
\end{equation}
\end{lemma}

\begin{proof}
For $p=2$, the claim  (\ref{eq:averaged-multiplier}) follows from Plancherel's theorem and Cauchy-Schwarz inequality, while for $p$ close to $1_+$ the claim follows (ignoring the bounded $|\Omega_m(\xi;\delta)|$-factor) from the assumption that $m(\xi,v)$
satisfies the truncation property uniformly in $v$ (and in fact, the end point $p=1$, with the usual ${\mathcal H}^1$-replacement of $L^1$, can be treated by a refined argument along the lines of \cite{DLM91}). The general case of $1<p<2$ follows by interpolation.
\end{proof}

\noindent
\begin{remark}
Clearly, if $m(\xi,\cdot)$ and $c(\xi)$ are $L^p$-multipliers  then so is their product, and in particular, if $m$ has the truncation property then (\ref{eq:averaged-multiplier}) applies for $\psi(m(\xi,v)/\delta)c(\xi)$. 
\end{remark}

\subsection{Averaging lemma for homogeneous symbols}\label{subsec:hom}

We will present several versions of the averaging lemma.  The later versions will supercede the former, but for
pedagogical reasons we will start with the simpler case of homogeneous symbols.  In this case it is convenient to
use polar co-ordinates $\xi = |\xi| \xi'$, where $\xi' \in S^{d-1}$ is defined for all non-zero frequencies $\xi$ by
$\xi' := \xi/|\xi|$.

We begin with
\begin{avg}\label{thm:homogeneous-averaging}
Let $1\leq q\leq 2$ and let $g\in L^{q}_{\loc}(\Rvdx)$ if $q>1$ or let $g$ be a locally bounded measure,
$g\in {\cal M}(\Rvdx)$ if $q=1$.
Let $\eta,\mym \geq 0$ and let $f\in \Wsp{\Rvdx}, \ \sigma\geq0, \ 1 < p \leq 2$ solves the equation

\begin{equation}\label{eq:transport-again}
\op(\vnab) f(\vx) = \Lambda_x^\eta\partial_v^\mym g(\vx) \ {\rm in} \ {\dd}'(\Rvdx).
\end{equation}
Here, $\op(\nabla_x,\cdot)$ is a  differential operator with sufficiently smooth coefficients, $\op(\cdot,v)\in C^{\mym,\varepsilon>0}$, and let $\op(\vixi)$ be the corresponding symbol. We assume that $\op(\ixi, v)$ is homogeneous in $\xi$ of order $k, k>\sigma+\eta$, that the modified symbol $\op(\ixi',v)$ obeys the truncation property uniformly in $v$, and that it is nondegenerate in the sense that there exists an $\alpha, \ 0 <\alpha < (N+1)q'$, such that  
\begin{equation}\label{eq:non-degeneracy}
\sup_{\xi \in \Rd: |\xi| =1} | \myO | \lesssim \delta^\alpha,
\qquad \Omega_{\op}(\xi;\delta):= \Big\{v \in I \ : \ |\op(\vixi)| \leq \delta \Big\}.
\end{equation}

\noindent
Then, there exist $\theta=\theta_\alpha \in (0,1)$ and $s_\alpha=s(\theta_\alpha)> \sigma$ such that for all bump functions $\phi\in C_0^\infty(I)$, the averages $\overline{f}(x) :=\int f(\vx)\phi(v)dv$ belong to the Sobolev space 
$\Wsr{\Rdx}$ for all $s \in (\sigma,s_\alpha)$ and the following estimate holds
\begin{equation}\label{eq:regularity-estimate}
\|\overline{f}\|_{\Wsr{\Rdx}} \lesssim \|g\|^\theta_{\Lqloc(\Rvdx)}\cdot
\|f\|^{1-\theta}_{\Wsp{\Rvdx}}, \quad s \in (\sigma, s_\alpha), \ s_\alpha:= (1-\theta_\alpha)\sigma+\theta_\alpha(k-\eta).
\end{equation}
Here $\theta \equiv \theta_\alpha(p,q,N)$ and $r$ are given by
\begin{equation}\label{eq:theta}
\theta=\frac{\alpha/p'}{\alpha(1/p'-1/q')+\mym+1}, \quad \frac{1}{r}=\frac{1-\theta}{p}+\frac{\theta}{q}, \qquad 0< \theta<1.
\end{equation}
\end{avg} 

\begin{remark} It would be more natural to assume that the symbol $\op(\ixi,v)$ itself, rather than the modified symbol
$\op(\ixi',v)$, obeyed the truncation property, as it is typically easier to verify the truncation property for the unmodified symbol.  Indeed when we turn to more advanced versions of the averaging lemma (which rely on Littlewood-Paley theory, and which do
not assume homogeneity) we will work with the truncation property for the unmodified symbol.  However we choose to work here with
the modified symbol as it simplifies the argument slightly.
\end{remark}

\begin{proof}
We start with a smooth  partition of unity, $1\equiv \sum \psi_j(2^{-j}z)$ such that $\psi_0$ is a bump function supported inside the  disc $|z| \leq 2$ and   the other $\psi_j$'s are bump functions supported on the annulus $1/2<|z|<2$ (we note in passing  that the other $\psi_j$'s can be taken to be equal, so the index $j$ merely serves to signal their `action' on the shells, $2^{j-1}\leq|z|\leq 2^{j+1}$). We set
\bdm
f_j(x,v):= {\ff}_x^{-1} \psi_j\Big(\frac{\symb}{2^j\delta}\Big)\widehat{f}(\vxi), \qquad j=0,1,2 \ldots,
\edm
recalling that  $\xi':={\xi}/{|\xi|}$, and we consider the corresponding decomposition  $f= f_0+\sum_{j\geq1}f_j$.
We distinguish between two pieces, $f= f^{(0)}+f^{(1)}$ where $f^{(0)}:=f_0$ and $f^{(1)}:= \sum_{j\geq1}f_j$.
Observe that the $v$-support of $\widehat{f^{(0)}}$ is restricted to the degenerate set $\Omega_\op(\xi;\delta)$ whereas $\widehat{f^{(1)}}=\sum_{j\geq1}\widehat{f_j}$ offers a decomposition of the non-degenerate complement, $\Omega^c_\op(\xi;\delta)$. The free parameter $\delta$ is to be chosen later.
 
We start by noting that $f^{(0)}=f_0$ is associated with the  multiplier
$\psi_0(\symb/\delta)$.  Since $\op(\vixip)$ satisfies the truncation property, we 
can use Lemma \ref{basic-lp} and non-degeneracy assumption  (\ref{eq:non-degeneracy}), to obtain
\begin{equation}\label{eq:size-of-f1}
\|\overline{f^{(0)}}\|_{\Wsp{\Rdx}} \lesssim \supxi\Big|\myO \Big|^{1/p'} \cdot \|f\|_{\Wsp{\Rvdx}} \lesssim \delta^{\alpha/p'}\|f\|_{\Wsp{\Rvdx}}.
\end{equation} 
We turn to the other averages, $\overline{f_j}, \ j\geq1$ which make $f^{(1)}$. Since $\op(\ixi,\cdot)$ is homogeneous of order $k$, equation (\ref{eq:transport-again}) states that 
\bdm
\widehat{\Lambda_x^{k-\eta}\widehat f}(\vxi)=|\xi|^{k-\eta}\frac{|\xi|^\eta}{\op(\vixi)}\partial^\mym_v\widehat{g}(\vxi)= \frac{1}{\symb}\partial^\mym_v\widehat{g}(\vxi), 
\edm
and thus,
\bdm
\Lambda_x^{k-\eta}\overline{f_j} = \frac{1}{2^j\delta}{\ff}_x^{-1}\int_v \widetilde{\psi_j}\Big(\frac{\symb}{2^j\delta}\Big)
\partial^\mym_v\widehat{g}(\vxi)\phi(v)dv, \quad j=1,2,\ldots,
\edm
where $\widetilde{\psi_j}(z):=\psi_j(z)/z$ is a bump function much like $\psi_j$ is. Integration by parts then yields
\begin{eqnarray}\label{eq:2ndterm}
\ \ \Lambda_x^{k-\eta}\overline{f_j}   =  \frac{1}{(2^j\delta)^{N+1}}{\ff}_x^{-1}\int_v \widetilde{\psi_j}^{(N)}\Big(\frac{\symb}{2^j\delta}\Big)
\op_v^N(\vixip)&\!\!\! &\!\!\!\!\!\!\!\!\!\!\!\!\! \widehat{g}(\vxi)\phi(v)dv + \\
 &  + & \!\!\!\! \mbox{lower order or similar terms}. \nonumber
\end{eqnarray}
We can safely neglect the lower or similar order terms which involve (powers of) the bounded multipliers
$\partial^\ell_v\op(\vixip), \ell < N$ and we focus on the leading term in (\ref{eq:2ndterm}), associated  with the multipliers

\begin{equation}\label{eq:key-mult}
\widetilde{\psi_j}^{(N)}\Big(\frac{\symb}{2^j\delta}\Big)\op_v^N(\vixip).
\end{equation}
By H\"ormander-Mikhlin or Marcinkeiwicz multiplier theorems, $\op_v^N(\vixip)$ are bounded multipliers and hence (\ref{eq:2ndterm}) are upper-bounded by 
$$ \|\Lambda^{k-\eta}_x\overline{f_j}\|_{\Lqloc(\Rdx)}  \lesssim
\frac{1}{(2^{j}\delta)^{N+1}} \|\overline{\myM_j g}\|_{\Lqloc(\Rdx)}.$$
Here $\myM_j=\myM_{\widetilde{\psi_j}^{(N)}}$  are  the Fourier multipliers with symbol $\widetilde{\psi_j}^{(N)}({\symb}/{2^j\delta})$. 

We now fix $q>1$. By our assumption, $\symb$ satisfies the truncation property, and 
Lemma \ref{basic-lp} implies  $\|\overline{\myM_j g}\|_{\Lqloc(\Rdx)} \lesssim (2^{j+1}\delta)^{\alpha/q'}\|g\|_{\Lqloc(\Rdx)}$.
Adding together all the $f_j$'s, we find that $f^{(1)}=\sum_{j\geq 1} f_j$ satisfies

\begin{equation} \label{eq:size-of-f2}
\|\Lambda^{k-\eta}_x \overline{f^{(1)}}\|_{\Lqloc(\Rdx)}  \lesssim   \sum_{j\geq 1} \frac{1}{(2^{j}\delta)^{N+1}} (2^{j+1}\delta)^{\frac{\alpha}{q'}}\|g\|_{\Lqloc(\Rvdx)} \lesssim \delta^{\frac{\alpha}{q'}-(\mym+1)}\|g\|_{\Lqloc(\Rvdx)}. 
\end{equation}
Thus, if we  fix $t>0$ and choose $\delta$ to equilibrate the bounds  in (\ref{eq:size-of-f1}) and (\ref{eq:size-of-f2}),  
\[
\delta^{\alpha(1/p'-1/q')+\mym+1} \sim t \|g(\vx)\|_{\Lqloc}/\|f(\vx)\|_{W^{\sigma,p}_{loc}},
\]
then this tells us that
\[
 \inf_{\overline{f^{(0)}}+\overline{f^{(1)}}=\overline{f}}\Big[\|\overline{f^{(0)}}\|_{\Wsp{\Rdx}}+
	t\|\overline{f^{(1)}}\|_{\dot{W}^{k-\eta,q}_{loc}(\Rdx)}\Big] \lesssim   t^\theta 
\cdot\|g(\vx)\|^\theta_{\Lqloc(\Rvdx)}\cdot\|f(\vx)\|^{1-\theta}_{\Wsp{\Rvdx}}, 
\]
and the desired $W^{s,r}_{loc}$-bound follows for $s <\stheta$ with $\theta$ given in (\ref{eq:theta}).
The remaining case of $q=1$ can be converted into the previous situation using Sobolev embedding. In this case, $g$ being a measure, it belongs to  $W^{-\epsilon,q_\epsilon}$ for all $(\epsilon,q_\epsilon)$ such that

\[
g\in W^{-\epsilon,q_\epsilon}, \qquad \frac{d+2}{q_\epsilon'} < \epsilon <1 < q_\epsilon <\frac{d+2}{d+1},
\]   
and hence (\ref{eq:regularity-estimate}) applies for $s<(1-\theta)\sigma+\theta(k-\eta-\epsilon)$ and $\theta=\theta_\alpha(p,q_\epsilon,N)$;
we then let $\epsilon$ approach $0_+$ so that $q_\epsilon$ approaches arbitrarily close to $1_+$ to recover (\ref{eq:regularity-estimate}) with $s_\alpha$
and $\theta_\alpha(p,1,N)$.
\end{proof}

\begin{rem} As an example consider a (possibly, pseudo-)differential operator $\op(\nabla_x,\cdot)$ of order $k$ and let $f(x,v)\in W^{\sigma,2}_{loc}$ such that $\op(\nabla_x,v)f\in W^{\sigma-k+1,2}_{\loc}$. Assume that 
$\op_v(\ixi,v)\neq 0$ so that the nondegeneracy condition (\ref{eq:non-degeneracy}) holds with $\alpha=1$. Application of the averaging lemma \ref{thm:homogeneous-averaging} with $p=q=2, N=0$ and $\eta=k-\sigma-1$ then yields  the gain of half a derivative, $\overline{f}(x)\in W^{s,2}_{loc}$ with $s< \sigma/2+(\sigma+1)/2=\sigma+1/2$, in agreement with \cite[Theorem 2.1]{GG92}. The main aspect here is going beyond the $L^2$-framework, while allowing for general and  possibly different orders of integrability, $(f,g)\in (W^{\sigma,p},L^q)$.  
\end{rem} 

\begin{rem} The limiting case of the interpolation estimate (\ref{eq:regularity-estimate}), $\theta=1, s=s_\alpha$, corresponds to
Besov regularity $\overline{f}\in B^{s,r}_{t=\infty}(\Rdx)$. This regularity can be worked out using  a more 
precise bookkeeping of the Littlewood-Paley blocks.
For the transport case, $k=1$, it was carried out first in \cite[Theorem 3]{DLM91}, improved in \cite{Be94} and a final refinement with a secondary index $t=p$ can be found in \cite{DVP00}. This limiting case is encountered in the particular situation when  $k=\eta$, so that the interval $(\sigma, s_\alpha)$ `survives' at $\theta=1$. Here, one cannot expect for regularizing effect, but there is a persistence of relative compactness of the  mapping $g(\vx) \mapsto \overline{f}(x)$, \cite{PS98}. 
\end{rem}

\bigskip
 To gain a better insight into the last averaging lemma, we focus our attention on  the case where $f$ is a $W^{\sigma,p}_{loc}$-solution of
\begin{equation}\label{eq:meq1}
\op(\vnab) f(\vx) = \partial_v g(\vx), \quad g\in \Lqloc(\Rvdx), f\in \Wsp{\Rvdx},
\end{equation}
corresponding to the special case $\eta=0$ and $N=1$ in (\ref{eq:transport-again}). This case will suffice to cover all the single-valued applications we have in mind for the discussion 
in Sections \ref{sec:hyp}, \ref{sec:par} and \ref{sec:ellip} without the burden of carrying out an excessive amount of indices. 
The averaging lemma (\ref{thm:homogeneous-averaging}) implies that $\overline{f}(x)$ has Sobolev regularity of order 
$s < (1-\theta)\sigma+\theta k$,
\begin{equation}\label{eq:regularity-1st}
\|\overline{f}\|_{\Wsr{\Rdx}} \lesssim \|g\|^\theta_{\Lqloc(\Rvdx)}\cdot
\|f\|^{1-\theta}_{\Lploc(\Rdx)}, \quad  \theta\equiv\theta_\alpha:=\frac{\alpha/p'}{\alpha(1/p'-1/q')+2}.
\end{equation}

\noindent
The last regularity statement can be improved. To this end, we revisit the dyadic multipliers in (\ref{eq:key-mult}),
$\widetilde{\psi_j}^{(N)}\big({\symb}/{2^j\delta}\big)\op_v^N(\vixip)$. The key observation is that  
 $\op_v(\vixip)$ acts only on the \emph{subset} of $v$'s --- those which belong to  $v\in \Omega_\op(\xip;2^{j+1}\delta)$. Linking the size of $\op_v(\vixip)$ to that of $\symb$, we arrive at the following improved averaging regularity lemma.

\begin{avg}\label{thm:improved-averaging}
Let $f\in \Wsp{\Rvdx}, \ \sigma \geq 0,  1 < p\leq 2$, solves the equation 
\begin{equation}\label{eq:im-transport1}
\op(\vnab) f(\vx) = \partial_v g(\vx), \qquad g\in \left\{\begin{array}{ll}
L^{q}(\Rvdx), & 1 < q\leq 2,\\ \smallskip
{\cal M}(\Rvdx), & q=1. \end{array}\right.
\end{equation}
Let  $\op(\ixi,v)$ be the corresponding symbol. We assume that $\op(\ixi,v)$ is homogeneous in $\ixi$ of order $k, k>\sigma$, that the modified symbol $\op(\ixi',v)$ satisfies the  truncation property uniformly in $v$, and that it is nondegenerate in the sense that there exists an $\alpha, \ 0 <\alpha < q'$, such that (\ref{eq:non-degeneracy}) holds. Moreover,  assume that
\begin{equation}\label{eq:im-link}
\exists \mu \in [0,1] \ \hbox{s.t.}  \ \sup_{|\xi|=1} \ \sup_{v \in \Omega_\op(\xi;\delta)} |\op_v(\vixi)| \lesssim \delta^{\mu},
\quad \Omega_\op(\xi;\delta):=\left\{v\in I\ :\ |\op(\ixi,v)|\leq \delta\right\}.
\end{equation}

\noindent
Then, there exist $\theta=\theta_\alpha \in (0,1)$ and $r$ given by 
\begin{equation}\label{eq:im-theta}
\theta:=\frac{\alpha/p'}{\alpha(1/p'-1/q')+2-\mu}, \quad \frac{1}{r}:=\frac{1-\theta}{p}+\frac{\theta}{q}, 
\end{equation}

\noindent
such that for all bump functions $\phi\in C_0^\infty(I)$, the averages $\overline{f}(x):=\int f(\vx)\phi(v)dv$ belong to the Sobolev space 
$\Wsr{\Rdx}$ for all $s \in (\sigma,s_\alpha), \ s_\alpha=(1-\theta_\alpha)\sigma+\theta_\alpha k$ and (\ref{eq:regularity-estimate}) holds. 
\end{avg}

\noindent
For the proof we revisit (\ref{eq:2ndterm}) with $\eta=0, N=1$,
\[
\ \ \Lambda_x^{k}\overline{f_j}   =  \frac{1}{(2^j\delta)^{2}}{\ff}_x^{-1}\int_v \partial_z \widetilde{\psi_j}\Big(\frac{\symb}{2^j\delta}\Big)
\op_v(\vixip) \widehat{g}(\vxi)\phi(v)dv.
\]

\noindent
Its bound in (\ref{eq:size-of-f2}) can now be improved by the extra factor of $\big(2^{j+1}\delta\big)^{\mu}$ which follows from (\ref{eq:im-link}), yielding
 \bdm
\|\Lambda^{k}_x \overline{f_{\J}}\|_{\Lqloc(\Rdx)}  \lesssim   \sum_{j\geq 1} \frac{1}{(2^{j}\delta)^{2}} (2^{j+1}\delta)^{\frac{\alpha}{q'}+\mu}\|g\|_{\Lqloc(\Rvdx)} \lesssim 
\delta^{\frac{\alpha}{q'}-(2-\mu)}\|g\|_{\Lqloc(\Rvdx)}, 
\edm
and we conclude by arguing along the lines of the averaging lemma \ref{thm:homogeneous-averaging}. $\Box$

\begin{rem}
In the generic case of a symbol $\op(\ixi,v)$ which is analytic in $\ixi$ uniformly in $v$, the `degeneracy of order $\alpha$' in (\ref{eq:non-degeneracy}) implies that $\sup\left\{|\op_v(\vixi)|: \ |\xi|=1, |\op(\ixi,v)|\leq \delta\right\} \lesssim \delta^{1-\alpha}$. Thus, (\ref{eq:im-link}) holds with $\mu=1-\alpha$ and in this case, the velocity averaging holds 
with Sobolev-regularity exponent
\[
\overline{f}(x)\in W^{s,r}_{loc}(\Rdx), \quad s<s_\alpha:=(1-\theta_\alpha)\sigma+\theta_\alpha k, \ 
\theta_\alpha:= \frac{\alpha/p'}{\alpha(1/p'+1/q)+1}.
\]
\end{rem}
\begin{rem} The statement of the averaging lemma \ref{thm:improved-averaging} can be extended for $f\in L^p(\Rvdx)$ and $g\in \Lqloc(\Rvdx)$ in the full range of $1< p < \infty$ and $1 \leq q \leq \infty$.
The dual claim to (\ref{eq:averaged-multiplier}), based on $(L^2, BMO)$ interpolation for $2\leq p < \infty$ reads
\begin{equation} \label{eq:averaged-multiplier-extended}
\|\overline{\myM_\psi f}(x)\|_{L^p(\Rdx)}  \lesssim   \sup_\xi|\Omega_m(\xi;\delta)|^{\frac{1}{p}}\cdot\|f\|_{L^p(\Rvdx)}, \ 2\leq p < \infty. 
\end{equation}
This yields the same $W^{s,r}$-regularity as before
\begin{equation}
\|\overline{f}\|_{\Wsr{\Rdx}} \lesssim \|g\|^\theta_{\Lqloc(\Rvdx)}\cdot
\|f\|^{1-\theta}_{\Wsp{\Rvdx}}, \quad s \in (\sigma, (1-\theta_\alpha)\sigma+\theta_\alpha k)
\end{equation}
with $\theta=\theta_\alpha(\overline{p},\overline{q})$ given by
\bdm
\theta:=\frac{\alpha/\overline{p}}{\alpha(1/\overline{p}-1/\overline{q})+ 2-\mu}, \quad \frac{1}{r}:=\frac{1-\theta}{\overline{p}}+\frac{\theta}{\overline{q}}, \qquad 
\overline{p}:=\max(p,p'), \ \overline{q}:=\max(q,q').
\edm
\end{rem}

\subsection{An averaging lemma for general symbols}\label{subsec:gen}

We now turn our attention to averages involving general, not necessarily homogeneous symbols; as such the polar co-ordinate representation $\xi = |\xi| \xi'$ is no longer useful and will be discarded.  We focus on equations of the form   $\op(\vnab) f(\vx) = \partial_v g(\vx)$, corresponding to $(\eta,\mym)=(0,1)$ in (\ref{eq:transport-again}). 

\begin{avg}\label{avg}  
Let $f\in \Wsp{\Rvdx}, \ \sigma \geq 0,  1 < p\leq 2$, solves the equation 
\begin{equation}\label{eq:transport2}
\op(\vnab) f(\vx) = \partial_v g(\vx), \qquad g\in \left\{\begin{array}{ll}
L^{q}(\Rvdx), & 1 < q\leq 2,\\ \smallskip
{\cal M}(\Rvdx), & q=1. \end{array}\right.
\end{equation}
Let  $\op(\ixi,\cdot)$ be the corresponding symbol of degree $\leq k$ with sufficiently smooth $v$-dependent coefficients and assume it obeys the truncation property.
Denote
$$ \myo := \sup_{\xi \in \Rd: |\xi| \sim \dyM} \big| \Omega_\op(\xi;\delta)\big|,
\qquad \Omega_\op(\xi;\delta):=\big\{v  : \ |\op(\vixi)|\leq \delta\big\}$$
and suppose  the following non-degeneracy condition holds
\begin{equation}\label{wiff}
\exists \alpha, \beta>0 \ \ \hbox{s.t.}  \quad \myo \lesssim \Big(\frac{\delta}{\dyM^{\beta}}\Big)^{\alpha}, \quad \forall \delta > 0, \ \dyM \gtrsim 1.
\end{equation}
Moreover, assume that
\begin{equation}\label{eq:link}
\exists \lambda\geq 0 \ \hbox{and} \ \mu\in [0,1] \ \ \hbox{s.t.} \ \ \sup_{|\xi|\sim J} \ \sup_{v \in \Omega_\op(\xi;\delta)} |\op_v(\vixi)| \lesssim J^{\beta\lambda}\delta^{\mu}.
\end{equation}

\noindent
Then, for all bump functions $\phi\in C_0^\infty(I)$, the average $\overline{f}(x):=\int f(\vx)\phi(v)dv$ belongs to 
the Sobolev space $\Wsr{\Rdx}$ for $s \in (\sigma, s_{\alpha,\beta})$ and the following estimate holds

\begin{equation}\label{eq:r-avg}
 \Big\| \int f(\vx) \phi(v)\ dv \Big\|_{\Wsr{\Rdx}}
\lesssim \| f(\vx) \|_{\Wsp{\Rvdx}} + \| g(\vx)\|_{\Lqloc(\Rvdx)}.
\end{equation}
Here, $s_{\alpha,\beta}:=(1-\theta_\alpha)\sigma + \theta_\alpha\beta(2-\mu- \lambda)$ where $\theta \equiv \theta_\alpha$ and $r$ are given by
\begin{equation}\label{eq:r-theta2}
\theta:=\frac{\alpha/p'}{\alpha(1/p'-1/q')+2-\mu}, \quad \frac{1}{r}:=\frac{1-\theta}{p}+\frac{\theta}{q}, \qquad 0<\theta<1.
\end{equation}
\end{avg}

\begin{rem} \label{rem:im}
How does the last averaging lemma compare with the previous ones? We note that since $\deg\op_v(\ixi,\cdot) \leq k$, then the additional assumption (\ref{eq:link}) always holds with $\mu=0$ and $\beta\lambda=k$. Hence, the averaging lemma \ref{avg} with just the non-degeneracy condition (\ref{wiff})  yields 
 the  regularity  $\overline{f}(x)\in \Wsr{\Rdx}$ of (the reduced) order  $s \in (\sigma, s_{\alpha,\beta})$, where 
\begin{equation}\label{eq:r-theta}
s_{\alpha,\beta}:=(1-\theta_\alpha)\sigma + \theta_\alpha(2\beta- k), \quad \theta_\alpha:=\frac{\alpha/p'}{\alpha(1/p'-1/q')+2}.
\end{equation}

\noindent
Now, if in  particular, $\op(\ixi,\cdot)$ is homogeneous  of degree $k$, then  
\[
\omega_{\op}(\xi;\delta)  = \omega_{\op}\Big(\xip; \frac{\delta}{|\xi|^k}\Big), \qquad \xip:=\frac{\xi}{|\xi|};
\]
this shows that if  the non-degeneracy condition (\ref{eq:non-degeneracy}) of averaging lemma \ref{thm:homogeneous-averaging} holds, $\myo \sim \big({\delta}/{\dyM^k}\big)^{\alpha}$, then it implies (\ref{wiff})  with $\beta=k$, and we recover the homogeneous averaging lemma \ref{thm:homogeneous-averaging} with $(\eta,N)=(0,1)$, namely, the averages $\bar{f}$ gain regularity of order $s < (1-\theta)\sigma + \theta(2\beta-k)= (1-\theta)\sigma+\theta k$.  The only difference is that now the truncation property is
assumed on the unmodified symbol $\op(\ixi,v)$ rather than the modified one $\op(\ixi',v)$.

As for the averaging lemma \ref{thm:improved-averaging}, we first note that in the generic case of an homogeneous $\op(\cdot,v)$, the additional assumption (\ref{eq:link}) holds with $\lambda=\alpha$ and $\mu=1-\alpha$,
\begin{equation}\label{eq:im-alpha}
\sup_{|\xi|\sim J} \ \ \sup_{\{v\in I:\ |\op(\vixi)|\leq \delta\} } |\op_v(\vixi)| \lesssim J^{\beta\alpha}\delta^{1-\alpha}.
\end{equation}
Indeed, all the homogeneous examples discussed in S3ections \ref{sec:hyp},\ref{sec:par} and \ref{sec:ellip} below, employ the averaging lemma \ref{thm:improved-averaging} with these parameters which yield $W^{s,r}$-regularity of order
$s<(1-\theta_\alpha)\sigma +\theta_\alpha\beta$,
\begin{equation}\label{eq:opts}
\overline{f}(x)\in W^{s,r}_{loc}(\Rdx), \quad s<(1-\theta_\alpha)\sigma +\theta_\alpha\beta, \ \ \theta_\alpha=\frac{\alpha/p'}{\alpha(1/p'+1/q)+1}
\end{equation}
In the particular case of  $\op$ being homogeneous of order $k$ then $\beta=k$ and we recover the  averaging lemma \ref{thm:improved-averaging} (except that the truncation hypothesis is now assumed on the unmodified symbol).
\end{rem}

\begin{proof}
We begin by noting that we can safely replace the non-degeneracy condition (\ref{wiff}) with a slightly weaker one, namely 
\begin{equation}\label{eq:betterwiff}
\exists \alpha, \beta>0, \  \hbox{and} \ \eps>0 \ \ \hbox{s.t.}  \quad \myo \lesssim \Big(\frac{\delta^{1+\eps}}{\dyM^{\beta}}\Big)^{\alpha}, \quad \forall \delta > 0, \dyM \gtrsim 1,
\end{equation}
and still retain the same gained of regularity of order $s <s_{\alpha,\beta}$.
This can be achieved by replacing the values of $\alpha$  in (\ref{wiff}) by  $\alpha/(1+\eps)$  and then absorbing $\eps$ into a slightly smaller  order of regularity, $(1-\eps)s_{\alpha,\beta}$ dictated by $\theta_{\alpha/(1+\eps)}$. The extra $\eps$-power of $\delta$ will be needed below to insure simple summability, which probably could be eliminated by a more refined argument involving Besov spaces, along the lines of \cite{DLM91}. 

Next, we break up $f$ into Littlewood-Paley pieces, 
$$f = f_0 + \sum_{\operatorname{dyadic} \ \!  \dyM'\!s \gtrsim 1 } f_{\dyM},$$ 
so that $\widehat f_{\dyM}(\vxi)$, the spatial Fourier transform  of $f_{\dyM}(\vx)$, is supported for frequencies $|\xi| \sim \dyM$, and $\hat{f}_0$  has  support in $|\xi| \lesssim 1$.
Since $f_0$ is a smooth average of $f$ at unit scales, the contribution of $f_0$ is easily seen to be acceptable.  By giving up an $\eps$ in the index $s_{\alpha,\beta}$  we may thus reduce (\ref{eq:r-avg}) to a single value of $\dyM$.  It thus suffices to show that
\[
\dyM^s \| \overline{f_{\dyM}} \|_{L^r_{loc}(\Rd)}
\lesssim \dyM^\sigma \| f_{\dyM} \|_{\Lploc(\Rvdx)} + \| g \|_{\Lqloc(\Rvdx)},
\quad \overline{f_{\dyM}}:=\int_v f_{\dyM}(\vx)\phi(v) dv,
\]
for each $\dyM \gtrsim 1$.

Fix $\dyM$.   Because of the local nature of Littlewood-Paley  projections when $\dyM \gtrsim 1$ we may replace the localized $L^p$ norms with global norms.  Actually we may replace $L^r$ by weak $L^r$ since we may pay another $\eps$ in the index $s_{\alpha,\beta}$  to improve this.  By the duality of weak $L^r$ and $L^{r',1}$ it thus suffices to show that
\begin{equation}\label{ching}
|\langle \overline{f_{\dyM}}, \chi_E \rangle| \lesssim \dyM^{-s} |E|^{1/r'} 
\end{equation}
for all sets $E \subset \Rd$ of finite measure, where we have normalized
$\| f_{\dyM} \|_{L^{p}(\Rvdx)} \lesssim \dyM^{-\sigma}$ and $\| g \|_{\Lqloc(\Rvdx)} \lesssim 1$.

\medskip\noindent
We now decompose the action in $v$-space of each of the Littlewood-Paley pieces (rather than a decomposition $f$ itself used in lemma \ref{thm:homogeneous-averaging}),
\[
 f_{\dyM}(\vx) = \sum_{\operatorname{dyadic} \ \!  \delta'\!s \lesssim \dyM^k } \psi\Big(\frac{\op(\vnab)}{\delta}\Big) f_{\dyM}(\vx), \qquad
\psi\Big(\frac{\op(\vnab)}{\delta}\Big): ={\ff}_x^{-1}\psi\Big(\frac{\op(\vixi)}{\delta}\Big){\ff}_x.
\]
Here, $\psi(z)$ is a bump function on $\C$ supported on the region $|z| \sim 1$.
It suffices to estimate 
\begin{equation}\label{wookie}
\Big\langle \int  \psi\Big(\frac{\op(\vnab)}{\delta}\Big) f_{\dyM}(\cdot,v) \phi(v) dv, \chi_E \Big\rangle
\end{equation}
with a summable decay as the dyadic $\delta \to 0$, so that (\ref{ching}) holds.

\noindent
By our assumption, $\op(\vixi)$ satisfies the truncation property uniformly in $v$, hence by Lemma \ref{basic-lp} we see that for all $1 < p \leq 2$,
\begin{equation}\label{av}
\Big\| \int  \psi\Big(\frac{\op(\vnab)}{\delta}\Big) f_{\dyM}(\vx) \phi(v)dv \Big\|_{\Lploc(\Rdx)}
\lesssim \myo^{1/p'} \| f_{\dyM} \|_{\Lploc(\Rvdx)}.
\end{equation}
From (\ref{av}) and H\"older we may thus estimate (\ref{wookie}) by
\begin{equation}\label{bound-1}
\Big|\Big\langle \int  \psi\Big(\frac{\op(\vnab)}{\delta}\Big) f_{\dyM}(\cdot,v) \phi(v) dv, \chi_E \Big\rangle\Big| \lesssim \dyM^{-\sigma} \myo^{1/p'}  |E|^{1/p'}.
\end{equation}
On the other hand, thanks to equation (\ref{eq:transport2}) we can write
$$ \psi\Big(\frac{\op(\vnab)}{\delta}\Big) f_{\dyM}(\vx) = \widetilde \psi\Big(\frac{\op(\vnab)}{\delta}\Big) \frac{1}{\delta} \frac{\partial}{\partial v} g_{\dyM}(\vx)$$
where $\widetilde \psi(z) := \psi(z)/z$ and the $g_J$'s are the corresponding Littlewood-Paley dyadic pieces of $g$. We thus have
\begin{equation}\label{eq:f2g}
 \int  \psi\Big(\frac{\op(\vnab)}{\delta}\Big) f_{\dyM}(\vx)\phi(v) dv = \frac{1}{\delta} \int  
\widetilde \psi\Big(\frac{\op(\vnab)}{\delta}\Big) \frac{\partial}{\partial v} g_{\dyM}(\vx)\phi(v) dv.
\end{equation}
We now integrate by parts to move the $\partial/\partial_v$ derivative somewhere else.  We will assume that the derivative hits $\widetilde \psi(\op(\vnab)/\delta)$, as the case when the derivative hits the bump function $\phi(v)$ is much better.
We are thus led to estimate 
\begin{equation}\label{eq:key}
\frac{1}{\delta^2} 
\Big| \Big\langle
\int \widetilde \psi_z\Big(\frac{\op(\vnab)}{\delta}\Big)  \op_v(\vnab) g_{\dyM}(x,v)\phi(v) dv, \chi_E \Big\rangle\Big|.
\end{equation}
Since $g_{\dyM}$ is localized to frequencies $\sim\dyM$, then by (\ref{eq:link}), the multiplier $\op_v(i\vxi)$ acts like a constant of order ${\mathcal O}(J^{\beta\lambda}\delta^\mu)$.  Also, $\widetilde \psi_z$ is a bump function much like $\psi$.  Thus we may modify (\ref{av})  --- with $p$ replaced by $q$ (assuming that $q>1$ and using the modified argument for the case $q=1$ as before), $\psi$ replaced by $\widetilde \psi_z$, and $f_{\dyM}$ replaced by $\op_v(\vnab) g_{\dyM}$, to estimate (\ref{eq:f2g})  by
\begin{equation}\label{bound-2}
\Big|\Big\langle \int  \psi\Big(\frac{\op(\vnab)}{\delta}\Big) f_{\dyM}(\cdot,v) \phi(v) dv, \chi_E \Big\rangle\Big| \lesssim \delta^{-(2-\mu)} \dyM^{\beta\lambda}  \myo^{1/q'} |E|^{1/q'}.
\end{equation}
Interpolating this bound with (\ref{bound-1}), we may bound (\ref{wookie}) by
$$
\delta^{-\theta(2-\mu)}  \dyM^{(1-\theta)(-\sigma) + \theta \beta\lambda} \myo^{1/r'} |E|^{1/r'}.$$
The parameterization in (\ref{eq:r-theta2}) dictates $\alpha/r'=\theta(2-\mu)$. Finally, we put the extra $\eps$ power into use: by (\ref{eq:betterwiff}), 
$\omega_\op(J;\delta)^{1/r'} \lesssim \big(\delta^{1+\eps}J^{-\beta}\big)^{\theta(2-\mu)}$ and hence the last quantity is bounded by 
\[
\delta^{\theta(2-\mu)\eps} \dyM^{- s_{\alpha,\beta} }|E|^{1/r'}, \qquad s_{\alpha,\beta}=(1-\theta)\sigma-\theta\beta(2-\mu- \lambda).
\]
Summing in $\delta$ and using the hypothesis that $s< s_{\alpha,\beta}$ we obtain (\ref{ching}).
\end{proof}

\subsection{Velocity averaging for first and second order symbols}\label{subsec:fssymbols}

To apply the velocity averaging \ref{avg} we need to find out which multipliers $m(\xi)$ have the truncation property.  Fortunately, there are a large classes of such multipliers.

First of all, it is clear that the multipliers $m(\xi) = \xi \cdot e_1$ and $m(\xi) = |\xi|^2$ have the truncation property, as in these cases the Fourier multipliers are just convolutions with finite measures.  Now, observe that if $m(\xi)$ has the truncation property, then so does $m(L(\xi))$ for any invertible linear transformation $L$ on $\Rd$, with a bound which is uniform in $L$.  This is because the $L^p$ multiplier class is invariant under linear transformations.

Because of this, we see that the multipliers
$$ m_1(\xi) = \bba(v) \cdot i\xi$$
and
$$ m_2(\xi) = \langle \bbb(v)\xi, \xi\rangle$$
have the truncation property uniformly in $v$, where $\bba(v)$ are arbitrary real coefficients, and $\bbb(v)$ is an arbitrary elliptic bilinear form with real coefficients. 

From the H\"ormander-Mikhlin or Marcinkeiwicz multiplier theorems and the linear transformation argument 
one can also show that $m_1(\xi')$ has the truncation property uniformly in $v$.  These arguments go back to the discussion of \cite{DLM91}. The situation with $m_2(\xi')$ is less clear,
but fortunately we will not need to verify that these second order modified multipliers obey the truncation property since our averaging lemmata also work with a truncation property hypothesis on the unmodified multiplier.

Now we observe that if $m_1(\xi)$, $m_2(\xi)$ are real multipliers with the truncation property, then the complex multiplier $m_1(\xi) + i m_2(\xi)$ also has the truncation property.  The basic observation is that one can use Fourier series to write any symbol of the form
$$ \psi\Big( \frac{m_1(\xi) + i m_2(\xi)}{\delta} \Big)$$
as
$$ \sum_{j,k \in \Z} \widehat \psi(j,k) \widetilde \psi_j\Big( \frac{m_1(\xi)}{\delta} \Big) \widetilde \psi_k\Big( \frac{m_2(\xi)}{\delta} \Big)$$
where $\widetilde \psi_j(x) := e^{2\pi i j x} \widetilde \psi$ and $\widetilde \psi$ is some bump function which equals 1 on the one-dimensional projections of the support of $\psi$.  Since the $C^\ell$ norm of $\widetilde \psi_j$ grows polynomially in $j$ and $\widehat \psi(j,k)$ decays rapidly in $j$, $k$ (if $\psi$ is sufficiently smooth), we are done since the product of two $L^p$ multipliers is still an $L^p$ multiplier.

\section{Nonlinear hyperbolic conservation laws}\label{sec:hyp}

Having developed our averaging lemmata, we now present some applications to nonlinear PDE.  We begin with the 
study (real-valued) solutions $\rho(\tx) = \rho(t,x_1, \ldots, x_d) \in L^\infty(\Rt \times \Rdx)$ of multidimensional scalar conservation laws

\begin{equation}\label{eq:cl}
\frac{\partial}{\partial t}\rho(t,x) + \sum_{j=1}^d\frac{\partial}{\partial x_j}A_j(\rho(t,x)) = 0 , \qquad \hbox{in} \ 
\dd'(\Rt\times \Rdx).
\end{equation}

\noindent
We abbreviate \eqref{eq:cl} as $\rho_t + \nabla_x\cdot\bbA(\rho)=0$ where $\bbA$ is the vector of $C^{2,\epsilon}$-spatial fluxes, $\bbA:=(A_1,A_2, \ldots,A_d)$. \newline
Let $\chi_\gamma(v)$ denote the velocity indicator function
\[
\chi_\gamma(v)=\left\{\begin{array}{cl} 1 & \hbox{ if} \ 0 < v \leq \gamma \\ -1 & \hbox{ if} \ \gamma\leq v < 0 \\0 & 
\hbox{ otherwise}. \end{array}\right.
\]
We say that $\rho(\tx)$ is a \emph{kinetic solution} of the conservation law (\ref{eq:cl}) if  the
corresponding distribution function, $\chi_{\rho(\tx)}(v)$, satisfies the transport equation
\begin{equation}\label{eq:kin-transport}
\partial_t \chi_{\rho(\tx)}(v) +\bba(v)\cdot \nabla_x \chi_{\rho(\tx)}(v) = \partial_v m(\txv) \quad \hbox{in} \ \dd'(\Rt\times\Rvdx),
\end{equation}
for some nonnegative measure, $m(\txv)\in {\cal M}^+(\Rt \times \Rvdx)$. Here, $\bba(v)$ is the vector of transport velocities, $\bba(v):=(a_1(v), \ldots, a_d(v))$ where $a_j(\cdot):=A'_j(\cdot), j=1,2, \ldots d$. 
The regularizing effect associated with the proper notion of  \emph{nonlinearity} of the conservation law (\ref{eq:cl}) 
was explored in \cite{LPT94a} through the averaging properties of an underlying kinetic formulation. 
For completeness, we include here a brief description which will serve our discussion on nonlinear parabolic and elliptic 
equations in the next sections and we refer to \cite{LPT94a} for a complete discussion. 

The starting point are the \emph{entropy inequalities} associated with (\ref{eq:cl}), 
\[
\partial_t \eta(\rho(\tx)) +\nabla_x \cdot \Aeta(\rho(\tx)) \leq 0  \ \ \hbox{in} \ \ \dd'(\Rt\times \Rdx).
\]
Here,  $\eta$ is an arbitrary \emph{entropy function} (a convex function from $\R$ to $\R$) and $\Aeta:=(A^\eta_1,\ldots,A^\eta_d)$ is the corresponding vector of   
entropy fluxes, $\Aetaj(\rho):= \int^\rho \eta'(s)A'_j(s)ds, \ j=1,2,\ldots d$.
A function $\rho\in L^\infty$ is an \emph{entropy solution} 
if it satisfies the entropy inequalities for \emph{all} pairs $(\eta,\Aeta)$
induced by convex entropies $\eta$. Entropy solutions are precisely those solutions which are \emph{realizable} as vanishing viscosity limit solutions and  
are uniquely determined by their $L^\infty\cutL1$-initial data, $\rho_0(x)$, prescribed at $t=0$, e.g., \cite{Lax73}. A decisive role is played by the one-parameter family of \emph{Kru\v{z}kov entropy pairs}, $(\eta(\rho;v), \Aeta(\rho;v))$,
parameterized by $v\in \R$,
\[
\eta(\rho;v):=|\rho-v|, \quad \Aetaj(\rho;v):=\sgn(\rho-v)(A_j(\rho)-A_j(v)).
\]
Kru\v{z}kov entropy pairs lead to a complete $L^1$-theory of existence, uniqueness and stability of first-order quasilinear conservation laws, \cite{Kr70}. 
We turn to the kinetic formulation. We define the distribution, $m(\txv)=m_{\rho(\tx)}(v)$ by the formula
\begin{equation}\label{eq:defg}
m(\txv):= -\left[\partial_t \frac{\eta(\rho;v)-\eta(0;v)}{2} +\nabla_x \cdot \Big(\frac{\Aeta(\rho;v)-\Aeta(0;v)}{2}\Big)\right].
\end{equation}  
The entropy inequalities tell us that the distribution $m=m_\rho$ is in fact a nonnegative measure, 
$m(\txv)\in {\cal M}^+(\Rt\times\Rvdx)$. Next, we differentiate (\ref{eq:defg}) with respect to $v$: a straightforward computation yields 
that $\chi_{\rho(\tx)}(v)$ satisfies the kinetic transport equation (\ref{eq:kin-transport}). This reveals the interplay between Kru\v{z}kov entropy inequalities and the underlying kinetic formulation; for nonlinear conservation laws, kinetic solutions coincide with the entropy solutions, \cite{PT91, LPT94a}. Observe that by velocity averaging we recover the macroscopic quantities associated with  the entropy solution $\rho$,
\[
\int_v \chi_{\rho}(v) \phi(v) dv = \Phi(\rho)
\]
where  $\Phi(\rho):=\int_{s=0}^\rho \phi(s) ds$ is the primitive of $\phi$.  In particular one
can recover $\rho$ itself by setting $\phi(v)=1_{[-M, M]}(v), \ M=\|\rho\|_{L^\infty}$.

We now use the averaging lemma (\ref{thm:improved-averaging}) to study the regularity of $\rho$. 
To this end we first extend (\ref{eq:kin-transport})
over the full $\Rtdx$-space, using a $C_0^\infty\Rt$-cut-off function,  $\psi\equiv 1$ for $t\geq \epsilon$, so that $f(\txv):=\chi_{\rho(\tx)}(v)\psi(t)$ and $g(\txv):=m_\rho(\txv)\psi(t)+\chi_{\rho(\tx)}(v)\partial_{t}\psi(t)$ satisfy 
\[
\partial_t f(\txv) +\bba(v)\cdot \nabla_x f(\txv) = \partial_v g(\txv) \ \ \hbox{in}\ \ \dd'(\Rtdxv), 
\quad   g\in {\cal M}(\Rtdxv).
\] 
Set $I:=[\inf \rho_0, \sup\rho_0]$ and assume that the first-order symbol is non-degenerate, (\ref{eq:non-degeneracy}), (\ref{eq:im-link}), namely, 
\begin{equation}\label{eq:1st-symb}
\exists \alpha \in (0,1) \ \  \hbox{s.t.} \  \sup_{{\tau}^2+|\xi|^2=1}|\Omega_\bba(\xi;\delta)|\lesssim \delta^\alpha, \qquad 
\Omega_\bba(\xi;\delta):=\Big\{v\in I :  |\tau+\bba(v)\cdot\xi|\leq \delta\Big\}, 
\end{equation}
and
\begin{equation}\label{eq:1st-symb-mu}
\exists  \mu\in [0,1] \ \  \hbox{s.t.} \ \ \sup_{|\xi|=1} \ \sup_{\Omega_\bba(\xi;\delta)} |\bba'(v)\cdot\xi| \lesssim \delta^{\mu}.
\end{equation} 
 
We apply the  averaging lemma  \ref{thm:improved-averaging} for first-order symbols, $k=1$, with $q=1, p=2$ and  $\sigma=0$, to find that
$\overline{f}(x)$ and hence $\rho(x)$ belong to $\Wsr{\Rtdx}$,
\[
\rho(t,x)\in \Wsr{(\epsilon,\infty)\times\Rdx}, \qquad s < \theta_\alpha:=\frac{\alpha}{\alpha+4-2\mu}, \ r:=\frac{\alpha+4-2\mu}{\alpha+2-\mu}.
\]
At this stage, we invoke the monotonicity property of  entropy solutions, which implies  that  for $s<1$, $\|\rho(t,\cdot)\|_{W^{s,1}_{loc}(\Rdx)}$ is nonincreasing, and we  deduce that $\rho(t,\cdot) \in W^{s,1}(\Rdx)$
for $t>\epsilon$. We conclude that the entropy solution operator associated with the nonlinear conservation law (\ref{eq:cl}),(\ref{eq:1st-symb}),  $\rho_0(\cdot) \mapsto \rho(t,\cdot)$, has a \emph{regularizing effect}, mapping $L^\infty(\Rdx)$ into $W^{s,1}_{loc}(\Rdx)$,
\[
\forall t\geq \epsilon >0: \quad \rho(t,\cdot) \in W^{s,1}_{loc}(\Rdx), \quad \ s < s_1, \ s_1:=\theta_\alpha=\frac{\alpha}{\alpha+4-2\mu}.
\]

Next, we use the bootstrap argument of \cite[\S3]{LPT94a} to deduce an improved regularizing effect.
The $W^{s,1}_{loc}(\Rtdx)$-regularity of $\rho(\tx)\psi(t)$ implies that $f(\txv)=\chi_{\rho(\tx)}(v)\psi(t)$ belongs to
$L^1(W^{s,1}(\Rtdx),\Rv)$; moreover, since $\partial_v\chi_\rho(v)$ is a bounded measure, $f \in L^1(W^{s,1}(\Rtdx),\Rv) \cap L^1(\Rtdx, W^{s,1}(\Rv))$ and hence
\[
f \in W^{s,1}_{loc}(\Rtdxv), \quad \forall s<s_1. 
\]
Interpolation with the obvious $L^\infty$-bound of $f$ then yields that $f\in W^{s,2}_{loc}(\Rtdxv)$ for all
$s <s_1/2$. Therefore, the averaging lemma \ref{thm:improved-averaging} applies to $f=\chi_{\rho(\tx)}(v)\psi(t)$ with $q=1, p=2$ and $\sigma=s_1/2$, implying that $\rho(t,\cdot)$ has improved  $W^{s,1}_{loc}$-regularity
 of order $s<s_2=(1-\theta_\alpha)s_1/2+\theta_\alpha$. Reiterating this argument yields the fixed point
$s_k\uparrow s_\infty=2\theta_\alpha/(1+\theta_\alpha)$ and we conclude with a regularizing effect
\[
\forall t\geq \epsilon >0: \quad  \rho_0\in L^\infty\cutL1(\Rdx) \mapsto \rho(t,\cdot) \in W^{s,1}_{loc}(\Rdx), \quad \ s < \frac{\alpha}{\alpha+2-\mu}.
\]

\noindent
As indicated earlier in remark \ref{rem:im}, in the generic case, $\mu=1-\alpha$,
\begin{equation}\label{eq:1st-symb-alpha}
 \sup_{|\xi|=1} \sup_{\{ v\in I: |\tau+\bba(v)\cdot\xi|\leq \delta\}} |\bba'(v)\cdot\xi| \lesssim \delta^{1-\alpha},
\end{equation} 
which yields $W^{s,1}$-regularizing effect of order $s<\alpha/(2\alpha+1)$. This improves the previous regularity result \cite[Theorem 4]{LPT94a} of order $s<\alpha/(\alpha+2)$, corresponding to $\mu=0$.  
\noindent
We can extend the last statement for general $L^p_{loc}$ initial data.
Recall that the entropy solution operator associated with (\ref{eq:cl}) is $L^1$-contractive. We now invoke a general \emph{nonlinear} interpolation argument of J.-L. Lions, e.g., \cite[Interpolation Theory, Lecture 8]{Ta00}; namely, if a possibly nonlinear $T$ is Lipschitz on $X$ with a Lipschitz constant $L_X$ and maps boundedly  $Y_1 \mapsto Y_2$ with a bound $B_Y$,  then one verifies that the corresponding $K$-functionals satisfy $K(Tx,t;X,Y_2)\leq L_XK(x,tB_Y/L_X;X,Y_1)$, and hence $T$ maps $[X,Y_1]_{\theta,q} \mapsto [X,Y_2]_{\theta,q}$. Consequently, the entropy solution operator maps $[L^1,L^\infty]_{\theta,q} \mapsto [L^1, W^{s,1}_{loc}]_{\theta,q}, \ 
0<\theta<1 < q$, and we conclude

\begin{cor}\label{cor:lpdata}
Consider the nonlinear conservation law (\ref{eq:cl}) subject to $L^p\cutL1$-initial data, $\rho(0,x)=\rho_0(x)$. Assume the non-degeneracy condition of order $\alpha$, (\ref{eq:1st-symb}), (\ref{eq:1st-symb-alpha}) holds over arbitrary finite intervals $I$. Then $\rho(\tx)$ gains a regularity of order $s/p'$, 
\[
\forall t\geq \epsilon >0: \quad  \rho_0\in L^p\cutL1(\Rdx) \mapsto \rho(t,\cdot) \in W^{s,1}_{loc}(\Rdx), \quad \ s < \frac{\alpha}{(2\alpha+1)p'}.
\]
\end{cor}

\noindent
The study of regularizing effects in one- and two-dimensional \emph{nonlinear} conservation laws has been studied by a variety of different  approaches; an incomplete list of references  includes \cite{Ol63},\cite{Ta79},\cite{Ta87},\cite{EE93},\cite{TRB05}.

\bigskip
We close this section with three examples. Let $\ell \geq 1$ and consider the one-dimensional conservation law
\begin{equation}\label{eq:burgers}
\frac{\partial}{\partial t} \rho(\tx)+ \frac{\partial}{\partial x}\Big\{\frac{1}{\ell+1}\rho^{\ell+1}(\tx)\Big\}=0, \quad 
\rho_0\in [-M,M].
\end{equation}
It satisfies the non-degeneracy condition (\ref{eq:1st-symb}) with $\alpha=1/\ell$, hence 
$\rho(t,\cdot)_{|t>\epsilon}\in W^{s,1}_{loc}$ with $s<\alpha/(2\alpha+1)=1/(\ell+2)$. It is well-known, however,  that the entropy solution operator of  the inviscid Burgers' equation corresponding to $\ell=1$, maps $L^\infty \mapsto BV$, \cite{Ol63}. This shows that the regularizing effect of order 
$\alpha/(2\alpha+1)$ stated in corollary \ref{cor:lpdata} is not sharp (although the averaging argument is! consult \cite{DLW05} following \cite{JP02}). Accordingly, it was conjectured in \cite{LPT94a}  that (\ref{eq:1st-symb}) yields a regularizing effect of order $\alpha$.\newline 

Next, let $\ell, m \geq 1$ and consider the two-dimensional conservation law
\[
\frac{\partial}{\partial t} \rho(\tx)+ \frac{\partial}{\partial x_1}\Big\{\frac{1}{\ell+1}\rho^{\ell+1}(\tx)\Big\}
+\frac{\partial}{\partial x_2}\Big\{\frac{1}{m+1}\rho^{m+1}(\tx)\Big\}=0, \quad 
\rho_0\in [-M,M].
\]
If $\ell\neq m$ then (\ref{eq:1st-symb}) is satisfied with $\alpha=\min\{\frac{1}{\ell},\frac{1}{m}\}$ and we conclude
$\rho(t\geq \epsilon,\cdot)\in W^{s}_{loc}(L^1)$ with $s< \min \{\frac{1}{\ell+2},\frac{1}{m+2}\}$. If $\ell=m$, however, then there is
\emph{no regularizing effect} since $\tip+v^\ell\xip_1+v^m\xip_2 \equiv 0$ for $\tau=0, \xi_1+\xi_2=0$; indeed, $\rho_0(x-y)$ are steady solutions which allow oscillations to persist along $x-y=const$.
Other cases can be worked out based on their polynomial degeneracy; for example,
\[
\frac{\partial}{\partial t} \rho(\tx)+ \frac{\partial}{\partial x_1}\sin(\rho(\tx))
+\frac{\partial}{\partial x_2}\Big\{\frac{1}{3}\rho^{3}(\tx)\Big\}=0, \quad 
\rho_0\in [-M,M],
\]
has a non-degeneracy of order $\alpha=1/4$, yielding $W^{s,1}_{loc}$-regularity of order $s<1/6$.

\section{Nonlinear degenerate parabolic equations}\label{sec:par}

We are concerned with second-order, possibly degenerate parabolic equations in conservative form

\begin{equation}\label{eq:par}
\frac{\partial}{\partial t}\rho(\tx) + \sum_{j=1}^d\frac{\partial}{\partial x_j}A_j(\rho(\tx))
-\sum_{j,k=1}^d \frac{\partial^2}{\partial x_j\partial x_k}B_{jk}(\rho(\tx)) = 0 , \quad \hbox{in} \ 
\dd'(\Rt\times \Rdx).
\end{equation}
We abbreviate, $\rho_t + \nabla_x\cdot\bbA(\rho) + \trace\big(\nabla_x\otimes\nabla_x \ \bbB(\rho)\big)=0$ where $\bbB$ is the matrix $\bbB:=\big\{B_{jk}\big\}_{j,k=1}^d$. By degenerate parabolicity we mean  that the matrix $\bbB'(\cdot)$ is non-negative, $\big\langle \bbB'(\cdot)\xi,\xi\big\rangle \geq0, \ \forall \xi \in \R^d$. Our starting point are the
\emph{entropy inequalities} associated with (\ref{eq:par}),  such that for \emph{all} convex $\eta$'s,
\begin{equation}\label{eq:entropypar}
\partial_t \eta(\rho(\tx)) + \nabla_x\cdot \Aeta(\rho(\tx)) -\trace\Big(\nabla_x \otimes \nabla_x \ \Beta(\rho(\tx))\Big) \leq  0 \ \ \hbox{in} \ \ \dd'((0,\infty) \times \Rdx).
\end{equation}

\noindent
Here, $\Aeta$ is the same vector of hyperbolic entropy fluxes we had before, $\Aeta=(A^\eta_1,\ldots, A^\eta_d)$ and $\Beta$ is the matrix of parabolic entropy fluxes, $\Beta:=(B^\eta_{jk})_{j,k=1}^d, \  \Betajk(\rho):=\int^\rho \eta'(s) B_{jk}'(s)ds$. 
We turn to the kinetic formulation. Utilizing the Kru\v{z}kov entropies, $\eta(\rho;v):=|\rho-v|$, we define the distribution, $m(\txv)=m_{\rho(\tx)}(v)$,
\begin{eqnarray}
m(\txv)&:=& -\left[\partial_t \frac{\eta(\rho;v)-\eta(0;v)}{2} + 
\nabla_x \cdot \Big(\frac{\Aeta(\rho;v)-\Aeta(0;v)}{2}\Big)\right] + \nonumber \\
  & &  \quad + \ \trace\Big(  \nabla_x \otimes \nabla_x \frac{\Beta(\rho;v)-\Beta(0;v)}{2}\Big).
\end{eqnarray} 
The entropy inequalities tell us that $m(\txv)\in {\cal M}^+(\Rt\times\Rvdx)$ and differentiation with respect to $v$
yields the kinetic formulation,

\begin{equation}\label{eq:kinpar}
\partial_t \chi_{\rho(\tx)}(v) + \bba(v)\cdot\nabla_x \chi_{\rho(\tx)}(v) -\nabla_x^\top\cdot\bbb(v)\nabla_x \chi_{\rho(\tx)}(v) = \partial_v m(\txv), 
\end{equation}
for some nonnegative  $m\in{\cal M}^+$ which measures entropy+dissipation production. Here, $\bba$ is the same vector of velocities we had before, $\bba=\bbA'$, and $\bbb$ is the non-negative diffusion matrix, $\bbb:=\bbB'\geq 0$. 
The representation $\eta(\rho)-\eta(0)= \int \eta'(s)\chi_\rho(s)ds$ shows  that the kinetic formulation (\ref{eq:kinpar}) is in fact the equivalent dual statement of the entropy inequalities (\ref{eq:entropypar}). But neither of these statements settles the question of uniqueness, except for certain special cases, such as the isotropic diffusion, $B_{jk}(\rho)=B(\rho)\delta_{jk}\geq 0$, e.g., \cite{Ca99}, or special cases with mild singularities, e.g., a porous- media type one-point degeneracy, \cite{DiB93, Ta97}. 
The extension of Kru\v{z}kov theory to the present context of general parabolic equations with possibly non-isotropic diffusion was completed only recently in \cite{CP03}, after the pioneering work \cite{VH69}. Observe that the entropy production measure, $m$, consists of contributions from the hyperbolic entropy dissipation  and the parabolic dissipation of the equation, $m=m_\bbA+m_\bbB$.
The  solutions sought by Chen and Perthame in \cite{CP03}, $\rho\in L^\infty$,  require that their corresponding distribution function  $\chi_\rho$ satisfies (\ref{eq:kinpar})  with a \emph{restricted} form of parabolic defect measure $m_\bbB$: the restriction  imposed on $m_\bbB$  reflects a certain renormalization property of the mixed derivatives of $\rho$ (or more precisely, the primitive of $\sqrt{\bbb(\rho)}$).  Accordingly, we can refer to these Chen-Perthame solutions as \emph{renormalized solutions} with a kinetic formulation 
(\ref{eq:kinpar}). These renormalized kinetic solutions admit an equivalent interpretation as entropy solutions, \cite{CP03} and as dissipative solutions, \cite{PS05}.  A general $L^1$-theory of existence, uniqueness and stability can be found in  \cite{CK04}. For a recent overview with a more complete list of references on such convection-diffusion equations in divergence form we refer to \cite{Ch03}. The regularizing effect of such equations, however, is less understood. In \cite[\S5]{LPT94a} we used the kinetic formulation (\ref{eq:kinpar}) to prove that the solution operator, $\rho_0 \mapsto \rho(t,\cdot)$ is relatively compact under a generic non-degeneracy condition 
\[
\supxi |\Omega_\op(\xi;0)|=0, \quad \Omega_\op(\xi;0):=\left\{v \ : \ \tau+\bba(v)\cdot\xi =0, \ \langle\bbb(v)\xi,\xi\rangle =0 \right\}.
\]
A general compactness result in this direction can be can be found  \cite{Ge90}.
We turn to \emph{quantify}  the  regularizing effect associated with such kinetic solutions. We emphasize that our regularity results are based on the `generic' kinetic formulation (\ref{eq:kinpar}), but otherwise, they  are \emph{independent} of the additional information on the renormalized Chen-Perthame solutions encoded in their entropy production measure $m$. The extra restrictions of the latter will likely to yield even better regularity results than those stated below. We divide our discussion into two stages, in order to highlight different aspects of degenerate diffusion, in Section \ref{subsec:deg}, and the coupling with nonlinear convection, in Section \ref{subsec:cd}.

\subsection{Non-isotropic degenerate diffusion}\label{subsec:deg}
We consider the parabolic equation

\begin{equation}\label{eq:noniso}
\frac{\partial}{\partial t}\rho(\tx) -\sum_{j,k=1}^d \frac{\partial^2}{\partial x_j\partial x_k}B_{jk}(\rho(\tx)) = 0 , \quad \hbox{in} \ 
\dd'(\Rt\times \Rdx).
\end{equation}
Here, we ignore the hyperbolic part and focus on the effect of non-isotropic diffusion. The corresponding kinetic formulation
(\ref{eq:kinpar})  extended to the full $\Rtdx\times\Rv$ reads
\[
\partial_t f(\txv) - \nabla_x^\top \cdot\bbb(v)\nabla_x f(\txv) = \partial_v g(\txv), \quad
 f:=\chi_\rho\psi(t), \ g\in {\cal M}^+(\Rtdxv).
\]
Set $I:=[\inf\rho_0,\sup\rho_0]$. The corresponding symbol is  $\op(\tau,\ixi,v)=i\tau+\langle \bbb(v)\xi, \xi\rangle$ and it suffices to make the non-degeneracy assumption (\ref{eq:non-degeneracy}), on the second-order homogeneous part of the symbol $\op(0,\ixi,v)=\langle\bbb(v)\xi,\xi\rangle$. We make

\begin{equation}\label{eq:2nd-symb}
\exists \alpha \in (0,1) \ \ \hbox{s.t.} \ \sup_{|\xi|=1}\big|\Omega_{\bbb}(\xi;\delta)\big| \lesssim \delta^\alpha, \qquad 
 \Omega_{\bbb}(\xi;\delta):=\big\{v\in I: \ 0\leq \langle \bbb(v)\xi,\xi\rangle \leq \delta\big\}.
\end{equation}
and
\begin{equation}
\exists \mu \in [0,1] \ \ \hbox{s.t.} \ \ \sup_{|\xi|=1}\  \sup_{\Omega_{\bbb}(\xi;\delta)}\big|\langle\bbb'(v)\xi,\xi\rangle\big| \lesssim \delta^{\mu}
\end{equation}

\noindent
We apply the averaging  result \ref{avg} with $q=1,p=2, \sigma=0$ and $k=2$, to find  that $\overline{f}(\tx)$ and hence $\rho(\tx)$ belong to $\Wsr{\Rtdx}$,
\[
\rho(t,\cdot)\in \Wsr{(\epsilon,\infty)\times\Rdx}, \qquad s < 2\theta_\alpha, \ \theta_\alpha:=\frac{\alpha}{\alpha+4-2\mu}, \ r:=\frac{\alpha+4-2\mu}{\alpha+2-\mu}.
\]
We follow the hyperbolic arguments. The kinetic solution operator associated with (\ref{eq:par}) is $L^1$-contractive, hence 
$\|\rho(t,\cdot)\|_{W^{s,1}_{loc}(\Rdx)}$ is nonincreasing and we conclude that $\forall t>\eps$, $\rho(t,\cdot)$ has  $W^{s,1}_{loc}$-regularity of order  $s < s_1, \ s_1:=2\theta_\alpha=2\alpha/(\alpha+4-2\mu)$.
We then bootstrap. Since $\chi_\rho(v)\psi(t) \in W^{s,2}_{loc}(\Rtdxv)$ for all $s <s_1/2$, we can apply the averaging lemma \ref{thm:homogeneous-averaging}
with $\sigma=s_1/2$ leading to $W^{s,1}_{loc}$-regularity of order $s_2:=(1-\theta_\alpha)s_1/2+2\theta_\alpha$ with fixed point $s_k\uparrow s_\infty=4\theta_\alpha/(1+\theta_\alpha)$,

\begin{equation}\label{eq:parreg}
\forall t\geq\epsilon >0: \quad \rho_0\in L^\infty\cutL1(\Rdx) \mapsto \rho(t,\cdot) \in W^{s,1}_{loc}(\Rdx), \quad s < \frac{2\alpha}{\alpha+2-\mu}.
\end{equation}

We distinguish between two different types of degenerate parabolicity, summarized in the following two corollaries.
 
\begin{cor}\label{cor:parI}[Degenerate parabolicity I. The case of a full rank]
Consider the degenerate parabolic equation (\ref{eq:par}) subject to $L^\infty\cutL1$-initial data, $\rho(0,x)=\rho_0$.
Let $\lambda_1(v) \geq \lambda_2(v) \geq \ldots \lambda_d(v)\geq 0$ be the eigenvalues of $\bbb(v)$ and assume that 
$\lambda_d(v)\equiv\!\!\!\!\!\!\!/ \ 0$ over $I=[\inf\rho_0,\sup\rho_0]$. Then, $\big\langle \bbb(v)\xi,\xi\big\rangle \geq \lambda_d(v)|\xi|^2$ and $\rho(\tx)$ has a regularizing  of order $s <2\alpha/(\alpha+2-\mu)p'$, i.e., (\ref{eq:parreg}) holds 
with $\alpha$ and $\mu$ dictated by  $\lambda_d(v)\equiv \lambda_d(\bbb(v))$, 
\[
 |\Omega_\lambda(\delta)| \lesssim \delta^\alpha \ \ \hbox{and} \ \
\sup_{|\xi|=1}\sup_{v\in \Omega_\lambda(\delta)}|\langle\bbb'(v)\xi,\xi\rangle| \lesssim \delta^\mu, \qquad \Omega_\lambda(\delta):=\big\{v \in I\ : \ 0\leq \lambda_d(v) \leq \delta\big\}.
\]
\end{cor}

\noindent
Corollary \ref{cor:parI} applies to the special case of \emph{isotropic} diffusion,

\begin{equation}\label{eq:pariso}
\frac{\partial}{\partial t} \rho(\tx)- \Delta B(\rho(\tx))= 0, \quad B'(v)\geq0,
\end{equation}
subject to $L^\infty\cutL1$-initial data, $\rho(0,\cdot)=\rho_0$.  If $b(\cdot):=B'(v)$ is  degenerate of order $\alpha$ in the sense that,
\[
\big|\Omega_b(\delta)\big| \lesssim \delta^\alpha, \ \ \hbox{and} \ \ 
\sup_{v\in \Omega_b(\delta)} |b'(v)| \lesssim \delta^{1-\alpha}, \quad \Omega_b(\delta):=\big\{v \in I\ : \ 0\leq b(v) \leq \delta\big\},
\]
then corollary \ref{cor:parI} implies 
$\forall t\geq \epsilon: \rho(t,\cdot)\in W^{s,1}_{loc}, \ s <2\alpha/(2\alpha+1)$. For $L^p\cutL1$-data $\rho_0$, the corresponding solution $\rho(t,\cdot)$ gains  $W^{s,1}_{loc}$-regularity of order $s<2\alpha/(2\alpha+1)p'$ and we conjecture, in analogy with the hyperbolic case, that the non-degeneracy (\ref{eq:2nd-symb}) yields an improved  regularizing effect of order $2\alpha/p'$.
Existence, uniqueness and regularizing effects of the isotropic equation (\ref{eq:pariso}) were studied earlier in \cite{BC79, BC81a,BC81b}. The prototype is provided by the porous media equation, 
\begin{equation}\label{eq:pm}
\frac{\partial}{\partial t} \rho(\tx)-  \Delta\Big\{\frac{1}{n+1}|\rho^n(\tx)|\rho(\tx)\Big\}= 0, \quad \rho(0,x)=\rho_0(x)\geq 0, \ \rho_0\in L^\infty.
\end{equation}
The velocity averaging yields $W^{s,1}$-regularity of order $2/(n+2)$ and, as in the hyperbolic case, it does not recover the optimal  H\"{o}lder continuity in this case, e.g.  \cite{DiB93, Ta96}.
In fact, the kinetic arguments do not yield continuity.  Instead, our main contribution here is to the \emph{non-isotropic} case where  we conjecture the same gain of regularity  driven by $\lambda_d(\bbb(v))$, as the isotropic regularity 
driven  by $b(v)$.

We continue with the more subtle case where $\bbb(\cdot)$ does not have a full-rank, so that
\[
\exists \ell, \ 1\leq \ell <d \ : \ \lambda_1(v)\geq \ldots \lambda_\ell(v)\geq0, \quad \lambda_{\ell+1}(v)\equiv \ldots \equiv \lambda_d(v)\equiv0
\]
Despite this stronger degeneracy, there is still some regularity that can be `saved'. To demonstrate our point, we consider the 2D case.
\begin{cor}\label{cor:parII}[Degenerate parabolicity II. The case of a partial rank]
We consider the 2D degenerate equation 
\begin{equation}\label{eq:par2D}
\frac{\partial}{\partial t}\rho(\tx) - 
\Big\{\frac{\partial^2}{\partial x^2_1}B_{11}(\rho(\tx)) + \frac{\partial^2}{\partial x_1 \partial x_2}B_{12}(\rho(\tx))
+ \frac{\partial^2}{\partial x^2_2}B_{22}(\rho(\tx))\Big\}  =0.
\end{equation}
subject to $L^\infty\cutL1$-initial data, $\rho(0,x)=\rho_0$. Assume strong degeneracy, $b^2_{12}(v)\equiv 4b_{11}(v)b_{22}(v)$, so that $\lambda_2(\bbb(v))\equiv 0, \forall  v\in I=[\inf \rho_0,\sup_0]$. In this case, $\langle\bbb(v)\xi,\xi\rangle =\Big(\sqrt{b_{11}}(v)\xi_1 +\sqrt{b_{22}}(v)\xi_2\Big)^2$
and $\rho(t,\cdot)$ admits a $W^{s,1}_{loc}$-regularity of  of order $s<2\alpha/(\alpha+2-\mu)$ which is dictated by 
the non-degeneracy,
\begin{equation}\label{eq:parII-alpha}
\sup_{|\xi|=1}\Big|\Omega_{\bbb}(\xi;\delta)\Big| \lesssim \delta^\alpha, \qquad \Omega_{\bbb}(\xi;\delta):=\Big\{ v \in I : \ \big|\sqrt{b_{11}(v)}\xi_1+\sqrt{b_{22}(v)}\xi_2\big|^2 \leq \delta \Big\}, 
\end{equation}
and
\begin{equation}\label{eq:parII-mu}
\sup_{|\xi|=1}\ \sup_{v\in\Omega_{\bbb}(\xi;\delta)}\Big| \frac{b_{11}'(v)}{\sqrt{b_{11}}(v)}\xi_1+\frac{b_{22}'(v)}{\sqrt{b_{22}}(v)}\xi_2\Big| \lesssim \delta^{\mu-1/2}.
\end{equation}
\end{cor}

\noindent
We distinguish between two extreme scenarios.

(i) If $|b_{11}(v)| \gg |b_{22}(v)| \ \forall v\in I$, then the regularizing effect (\ref{eq:parreg}) holds with $(\alpha,\mu)$ dictated by $b_{22}(v)$,
\[
 \big|\Omega_{b_{22}}(\delta)\big| \lesssim \delta^\alpha \ \ \hbox{and} \ \
\sup_{v\in \Omega_{b_{22}}(\delta)}\big|b_{22}'(v)\big| \lesssim \delta^\mu, \qquad 
\Omega_{b_{22}}(\delta):=\big\{v \in I\ : \ 0\leq b_{22}(v) \leq \delta\big\}.
\]

(ii) If $b_{11}(v) \equiv b_{22}(v)\ \forall v\in I$, then there is \emph{no} regularizing effect since 
the symbol $\Big(\sqrt{b_{11}}(v)\xi_1 +\sqrt{b_{22}}(v)\xi_2\Big)^2$ vanishes for all $\xi_1\pm\xi_2=0$ (so that
(\ref{eq:parII-alpha}) is fulfilled with $\alpha=0$). Indeed, the equation (\ref{eq:par2D}), with $b_{11}(v)= b_{22}(v)=:B'(v)$, takes the form
\[
\frac{\partial}{\partial t} \rho(\tx)- \Big\{\frac{\partial^2}{\partial x^2_1} \pm 2\frac{\partial^2}{\partial x_1 \partial x_2} + \frac{\partial^2}{\partial x^2_2}\Big\}B(\rho(\tx))  =0,
\]
and we observe that $\rho_0(x\mp y)$ are steady solutions which allow for oscillations to persist along $x\mp y=const$.

\subsection{Convection-diffusion equations}\label{subsec:cd}
We begin with the one-dimensional case
\begin{equation}\label{eq:cdoned}
\frac{\partial}{\partial t}\rho(\tx)+ \frac{\partial}{\partial x}A(\rho(\tx)) -\frac{\partial^2}{\partial x^2}B(\rho(\tx))=0.
\end{equation}
We consider the prototype example of  high-order Burgers' type nonlinearity, $a(v):=v^\ell, \ell\geq 1$ combined with porous medium diffusion $b(v)=|v|^{n},  n\geq 1$. The corresponding symbol is given by $\op((\tau,\ixi),v)=i\tau+v^\ell i\xi + |v|^n\xi^2$. We study the regularity of this convection-diffusion equation using the averaging lemma \ref{avg}, which employs the size of the set
\[
\Omega_\op \big(J;\delta\big):=\left\{v \ \Big| \ J|\tau+v^\ell\xi|+ J^2|v|^n \xi^2 \leq \delta \right\}, \quad  \tau^2+\xi^2=1, \ J\gtrsim 1,  \ \delta \lesssim 1.
\]

\noindent
Comparing diffusion  vs. nonlinear convection effects, we can distinguish here between three different cases. Clearly, $\Omega_\op(J;\delta) \subset \Omega_b:= \left\{v : \ |v|^n \leq {\delta}/{J^2} \right\}$, hence
$\omega_\op(J;\delta)\lesssim (\delta/J^2)^{1/n}$ and  (\ref{wiff}) holds with $\alpha_b=1/n$ and $\beta_b=2$. We shall use this bound whenever $n\leq \ell$, which is the case dominated by the parabolic part of (\ref{eq:cdoned}). Indeed, in this case we have $(\delta/J)^{1/\ell} \gtrsim (\delta/J^2)^{1/n}$ which in turn yields  
\[
 \sup_{v\in\Omega_\op(J;\delta)}|\op_v\big((\tau,\ixi),v\big)| \lesssim 
\sup_{v\in \Omega_b} \big(J |v|^{\ell-1} + J^2  |v|^{n-1}\big) \lesssim J^{2/n} \delta^{1-1/n}.
\]
This shows that (\ref{eq:link})  holds with $(\lambda_b, \mu_b)=(\alpha_b, 1-\alpha_b)$  and velocity averaging implies $W^{s,1}_{loc}$-regularizing effect with Sobolev exponent $s < {\beta_b\alpha_b}/{(3\alpha_b+2)}$,
\[
s<s_n= \frac{2}{2n+3}.
\]
We also have, $\Omega_\op(J;\delta) \subset \Omega_a:= \{v: \ |v|^\ell\lesssim \delta/J\}$, so that
$\omega_\op(J;\delta)\lesssim (\delta/J)^{1/\ell}$, i.e., (\ref{wiff}) holds with $\alpha_a=1/\ell$ and $\beta_a=1$. We shall use this bound whenever $n\geq 2\ell$, which is the case driven by the hyperbolic part (\ref{eq:cdoned}). In this case, $(\delta/J)^{1/\ell} \lesssim (\delta/J^2)^{1/n}$, hence 
\[
  \sup_{v\in \Omega_\op(J;\delta)}|\op_v\big((\tau,\ixi),v\big)| \lesssim \sup_{v\in \Omega_a} \big(J |v|^{\ell-1} + J^2  |v|^{n-1}\big) \lesssim J^{1/\ell} \delta^{1-1/\ell},
\]
implying that (\ref{eq:link})  is fulfilled with $(\lambda_a,\mu_a) =(\alpha_a, 1-\alpha_a)$. The corresponding Sobolev exponent, $s< {\beta_a\alpha_a}/{(3\alpha_a+2)}$,  is then given by, 
\[
s< s_\ell:= \frac{1}{2\ell+3}.
\]
 Finally, for intermediate $n$'s, $\ell<n<2\ell$, we interpolate the previous two $\omega_\op$-bounds (which are valid for all $n$'s), 
\[
\omega_\op(J;\delta)\lesssim (\delta/J)^{(1-\zeta)/\ell}(\delta/J^2)^{\zeta/n}, 
\]
for some $\zeta\in [0,1]$, which we choose as $\zeta:=(n/\ell)-1$,  so that  (\ref{wiff}) holds with 
$\alpha=(1-\zeta)/\ell+\zeta/n, \beta\alpha=(1-\zeta)/\ell+2\zeta/n$ and 
(\ref{eq:link}) holds with $(\lambda,\mu)=(\alpha,1-\alpha)$. This then yields the Sobolev-regularity exponent, $s<\beta\alpha/(3\alpha+2)$, 
\[
s<\frac{n+(2\ell-n)\zeta}{3n+3(\ell-n)\zeta+2n\ell}, \quad 
\zeta:=(n/\ell)-1, \ \ell<n<2\ell. 
\]
An additional bootstrap argument improves this Sobolev exponent, $s<\beta\alpha/(2\alpha+1)$, and we summarize the three different cases in 

\noindent
\begin{cor}\label{cor:burgers-pm} The convection-diffusion equation
\begin{equation}\label{eq:reg-bpm}
\frac{\partial}{\partial t}\rho(\tx)+ \frac{\partial}{\partial x}\Big\{\frac{1}{\ell+1}\rho^{\ell+1}(\tx)\Big\} -\frac{\partial^2}{\partial x^2}\Big\{\frac{1}{n+1}|\rho^{n}(\tx)|\rho(\tx)\Big\}=0, \quad \rho_0\in [-M,M],
\end{equation} 
has a regularizing effect, $\rho_0\in L^\infty(\R_x) \mapsto \rho(t > \eps,\cdot)\in W^{s,1}_{loc}(\R_x)$,
of order $s<s_{\ell,n}$ given by
\[
s_{\ell,n}= \frac{n+(2\ell-n)\zeta_{\ell,n}}{2n+2(\ell-n)\zeta_{\ell,n}+n\ell}, \quad 
\zeta_{\ell,n}:=\left\{\begin{array}{ll} 0 & n\leq \ell \\ (n/\ell)-1 & \ell<n < 2\ell \\1 & n\geq 2\ell \end{array}\right. 
\]
\end{cor}
\noindent
We note that when $n\leq \ell$, then $\Omega_b \subset \Omega_a$ and (\ref{eq:reg-bpm}) is dominated by degenerate diffusion with a regularizing effect of order  $s_{\ell,n}=s_n=2/(n+2)$. Thus, we recover the same order of regularity  we met with the `purely diffusive' porous medium equation (\ref{eq:pm}). If  $ n \geq 2\ell$, however, then $\Omega_a\subset \Omega_b$ and it is the hyperbolic part  which dominates diffusion, driving the overall regularizing effect of (\ref{eq:reg-bpm}) with order $s_{\ell,n}=s_\ell=1/(\ell+2)$; we recover regularity with the same order which we  met with the `purely convective' hyper-Burgers' equation (\ref{eq:burgers}). Finally, in the intermediate `mixed cases', $\ell < n < 2\ell$, we find a regularity of a
(non-optimal) order 
\[
s_{\ell,n}= \frac{n\ell+(2\ell-n)(n-\ell)}{2n\ell-2(n-\ell)^2+n\ell^2}, \quad \ell < n < 2\ell.
\]

\bigskip
We turn to the multi-dimensional case (\ref{eq:par}). The regularizing effect
is determined by the size of the set
\[
\Omega\big( J;\delta)=\left\{v \ \Big| \ J|\tau+\bba(v)\cdot\xi|+ J^2\langle \bbb(v)\xi,\xi\rangle \leq \delta \right\}, \quad  \tau^2+|\xi|^2=1, \ J\gtrsim 1, \ \delta \lesssim 1.
\]

Assume that the degenerate parabolic part of the equation has a full-rank, so that the smallest eigenvalue of $\bbb(v),
\lambda(v)\equiv \lambda_d(\bbb(v))$ satisfies
\begin{equation}\label{eq:cdparreg}
|\Omega_{\bbb}(\delta)|\lesssim \delta^{\alpha_{\bbb}} \ \ \hbox{and} \ \sup_{|\xi|=1}\sup_{v\in \Omega_{\bbb}(\delta)} |\langle\bbb'(v)\xi,\xi\rangle| \lesssim \delta^{1-\alpha_{\bbb}}, \quad \Omega_{\bbb}(\delta):= \Big\{v \ \Big| \ 0\leq \lambda_d(\bbb(v))\leq \delta \Big\}.
\end{equation}  
In this case, $\myo \lesssim (\delta/J^2)^{\alpha_{\bbb}}$ which yields a gain of $W^{s,1}_{loc}$-regularity of order $s<2\alpha_{\bbb}/(2\alpha_{\bbb}+1)$.
If in addition, the hyperbolic part of the equation has a non-degeneracy of order $\alpha_{\bba}$, namely
\[
\sup_{\tau^2+|\xi|^2=1}\big|\Omega_{\bba}\big((\tau,\xi),\delta\big)\big|\lesssim \delta^{\alpha_{\bba}}  \quad \hbox{and} \quad   
  \sup_{|\xi|=1} \ \sup_{v\in \Omega_{\bba}(\delta)} |\bba'(v)\cdot\xi| \lesssim \delta^{1-\alpha_{\bba}}, 
\]
where $\Omega_{\bba}\big((\tau,\xi);\delta\big):= \Big\{v: \ |\tau+\bba(v)\cdot\xi| \leq \delta \Big\}$,
then we can argue along the lines of corollary \ref{cor:burgers-pm} to conclude that there is an overall $W^{s,1}_{loc}$-regularity of order dictated by the relative size of $2\alpha_{\bbb}/(2\alpha_{\bbb}+1)$ and $\alpha_{\bba}/(2\alpha_{\bba}+1)$. As an example, we have the following.

\begin{cor}
Consider the two-dimensional convection-diffusion equation 
\begin{equation}\label{eq:cdtwod}
\frac{\partial}{\partial t}\rho(\tx) +  \frac{\partial}{\partial x_1}\Big\{\frac{1}{\ell+1}\rho^{\ell+1}(\tx)\Big\} +
\frac{\partial}{\partial x_2}\Big\{\frac{1}{m+1}\rho^{m+1}(\tx)\Big\} 
 - \sum_{j,k=1}^2 \frac{\partial^2}{\partial x_j\partial x_k}B_{jk}(\rho(\tx)) =0,
\end{equation}
 with non-degenerate diffusion, $B'(v) \gtrsim |v|^n$. Then, its renormalized kinetic solution admits a $W^{s,1}_{loc}$-regularizing effect,  $\rho_0\in L^\infty \mapsto \rho(t,\cdot)\in W^{s,1}_{loc}(\R^2_x)$, of order 
\[
 s < \left\{\begin{array}{ll} 
s_{\ell,m}:={\displaystyle \min\Big(\frac{1}{\ell+2},\frac{1}{m+2}\Big)}, & \ \ \hbox{if} \ n \geq 2\max(\ell,m) \ \ \hbox{and} \ \ \ell\neq m, \\ 
s_n:={\displaystyle \frac{2}{n+2}}, & \ \ \hbox{if} \ n\leq \min(\ell,m) \ \ \hbox{or} \ \ \ell=m, \\ 
s_{\ell,m,n}\in [s_{\ell,m}, s_n], & \ \ \hbox{if} \ \min(\ell,m)< n < 2\max(\ell,m).
\end{array}\right.
\]
\end{cor}

\medskip
Finally, we close this section with a third example of a fully-degenerate equation
\begin{equation}\label{eq:cddtwod}
\frac{\partial}{\partial t}\rho(\tx) + \Big(\frac{\partial}{\partial x_1}
+\frac{\partial}{\partial x_2}\Big)A(\rho(\tx)) 
-\Big(\frac{\partial^2}{\partial x^2_1} -2\frac{\partial^2}{\partial x_1\partial x_2}+\frac{\partial^2}{\partial x^2_2}\Big)B(\rho(\tx))=0.
\end{equation}
In this case, there is a stronger, rank-one parabolic degeneracy with no regularizing effect from the purely diffusion part,
since $\langle\bbb(v)\xi,\xi\rangle\equiv 0, \ \forall \xi_1-\xi_2=0$, and no regularizing effect from the purely convection part where $\bba(v)\cdot \xi\equiv 0, \ \forall \xi_1+\xi_2=0$. 
Nevertheless, the combined \emph{convection-diffusion}  does have a regularizing effect as demonstrated in 

\begin{cor}\label{cor:mix}
Consider the two-dimensional convection-diffusion equation 
\begin{eqnarray}\label{eq:hypdifs}
\frac{\partial}{\partial t}\rho(\tx) & + & \Big(\frac{\partial}{\partial x_1}
+\frac{\partial}{\partial x_2}\Big)\Big\{\frac{1}{\ell+1}\rho^{\ell+1}(\tx) \Big\} \nonumber \\
& - &\Big(\frac{\partial^2}{\partial x^2_1} -2\frac{\partial^2}{\partial x_1\partial x_2}+\frac{\partial^2}{\partial x^2_2}\Big)\Big\{\frac{1}{n+1}|\rho^{n}(\tx)|\rho(\tx)\Big\}=0.
\end{eqnarray}
For $n\geq 2\ell$ it admits a regularizing effect,  $\rho_0\in L^\infty\cutL1 \mapsto \rho(t,\cdot)\in W^{s,1}_{loc}(\R^2_x)$, of order $s< 6/(2+2n-\ell)$. \end{cor}
\noindent
As before, the Sobolev exponent computed here is not necessarily sharp, and as a result  no gain of regularity is stated for $n < 2\ell$. 

For proof, we use the averaging lemma \ref{avg} with the usual $(p,q)=(2,1)$, which yields a Sobolev-regularity exponent of order
$s=\beta\alpha(2-\mu-\lambda)/(\alpha+4-2\mu)$, where $\alpha, \beta,\lambda$ and $\mu$ characterize the degeneracy of the  symbol associated with (\ref{eq:hypdifs}),
\[
\op\big((\tau,\ixi),v\big)=J i\tau+ J v^\ell i(\xi_1+\xi_2)+ J^2|v|^n|\xi_1-\xi_2|^2, \quad \tau^2+|\xi|^2=1.
\]
We consider first those $\xi$'s,  $|\xi|=1$ such that $|\xi_1-\xi_2|\leq 1/5$. Here we have (say)  $|\xi_1+\xi_2| \geq 1/10$, so that
$\Omega_\op\big(J;\delta\big) \subset \Omega_\bba=\big\{v: \ |v|^\ell \lesssim {\delta}/{J} \big\}$,
and (\ref{wiff}) holds with $(\alpha_a,\beta_a)=(1/\ell,1)$; moreover, the growth
of $\op_v$ sought in (\ref{eq:link}) is bounded by,
\[
\sup_{v\in \Omega_\op(J;\delta)}|\op_v\big((\tau,\ixi),v\big)| \lesssim J \Big(\frac{\delta}{J}\Big)^{(\ell-1)/\ell} + J^2  \Big(\frac{\delta}{J}\Big)^{(n-1)/\ell} \lesssim J^{\lambda} \delta^{\mu},
\]
where
\[ 
\left\{\begin{array}{llll} \lambda:={1}/{\ell} &  \hbox{and} &  \mu:= 1-{1}/{\ell}, & \ \ \hbox{if} \ n\geq 2\ell, \\ \lambda:=2-{(n-1)}/{\ell} & \hbox{and} &  {\mu:=(n-1)}/{\ell},
& \  \ \hbox{if} \ n <2\ell.\end{array}\right. 
\]
In particular, if  $n<2\ell$ then the Sobolev exponent vanishes since $2-\mu-\lambda=0$, and we cannot deduce any  regularizing effect in this case. If $n\geq 2\ell$, however,  we compute the Sobolev exponent as before,  $s=s_\ell=1/(2\ell+3)$.
 Next, we consider the case when  $|\xi_1-\xi_2|\geq 1/5$, so that  $\Omega_\op\big(J;\delta\big) \subset \Omega_\bbb=\big\{v: \ |v|^n \lesssim {\delta}/{J^2} \big\}$ and (\ref{wiff}) holds with $(\alpha_b,\beta_b)=(1/n,2)$. As for (\ref{eq:link}), we have
\[
\sup_{v\in \Omega_\op(J;\delta)}|\op_v\big((\tau,\ixi),v\big)| \lesssim J \Big(\frac{\delta}{J^2}\Big)^{(\ell-1)/n} + J^2  \Big(\frac{\delta}{J^2}\Big)^{(n-1)/n} \lesssim J^{2\lambda} \delta^{\mu},
\]
where
\[ 
\left\{\begin{array}{llll} \lambda=1/n & \hbox{and} &\mu= 1-1/n, & \ \hbox{if} \ n\leq \ell, \\ \lambda=1/2-(\ell-1)/n & \hbox{and} &  \mu=(\ell-1)/n, 
& \ \hbox{if} \ n >\ell.  \end{array}\right. 
\]
In particular, if $n\geq 2\ell$ we compute in this case a  smaller Sobolev exponent $s=3/(3-2\ell+4n)\leq s_\ell$. The regularity result follows from the  bootstrap argument we discussed earlier which yields the final Sobolev exponent
$s=2\beta\alpha(2-\mu-\lambda)/(\alpha+2-\mu)$. $\Box$

\section{Nonlinear degenerate elliptic equations}\label{sec:ellip}

We consider the nonlinear, possibly degenerate elliptic equation
\begin{equation}\label{eq:ell}
-\sum_{j,k=1}^d \frac{\partial^2}{\partial x_j\partial x_k}B_{jk}(\rho(x)) = S(\rho(x)) , \quad \hbox{in} \ 
\dd'(\Gamma), \quad \bbb(\cdot):=\bbB'(\cdot) \geq 0,
\end{equation}
augmented with proper boundary conditions along the $C^{1,1}$-boundary  $\partial\Gamma$. We assume that the nonlinear source term, $S(\rho)$, is further restricted so that blow-up is avoided. 

\noindent
We begin with formal manipulations, multiplying (\ref{eq:ell}) against $\eta'(\rho)$, and `differentiating by parts' to find 

\begin{eqnarray}
\lefteqn{-\sum_{j,k=1}^d \eta'(\rho) \frac{\partial^2}{\partial x_j\partial x_k}B_{jk}(\rho) -\eta'(\rho) S(\rho) =}\nonumber \\
& & = \lefteqn{-\sum_{j=1}^d \frac{\partial}{\partial x_j}\Big(\eta'(\rho) b_{jk}(\rho)\frac{\partial \rho}{\partial x_k}\Big)
+ \sum_{j,k=1}^d  \eta''(\rho)b_{jk}(\rho) \frac{\partial \rho}{\partial x_j}\frac{\partial \rho}{\partial x_k}-\eta'(\rho) S(\rho) =} \nonumber \\
& & = -\trace\Big(\nabla_x\otimes\nabla_x \ \Beta(\rho)\Big) + \eta''(\rho)\big\langle \bbb(\rho)\nabla_x\rho,\nabla_x\rho\big\rangle 
-\eta'(\rho) S(\rho). \nonumber
\end{eqnarray}

\noindent
We arrive at the \emph{entropy inequalities} associated with (\ref{eq:ell}), stating that sufficiently smooth solutions of (\ref{eq:ell}) satisfy, for all convex $\eta$'s,
\begin{equation}\label{eq:ellent}
\ \ -\trace\Big(\nabla_x\otimes\nabla_x \ \Beta(\rho)\Big) \leq  
\eta'(\rho) S(\rho), \quad {\Beta(\rho)}= \Big\{\Betajk(\rho):=\int_0^\rho\eta'(s)b_{jk}(s)ds\Big\}_{j,k=1}^d.
\end{equation}
 
So far, we have not  specified the notion of solutions for (\ref{eq:ell}) since it seems that relatively little is known about a general stability theory for degenerate equations such as (\ref{eq:ell}).
The difficulty lies with the type of degeneracy  which does not lend itself to standard elliptic regularity theory, because the $B_{jk}$'s degenerate dependence on $\rho$,  nor does it admit the regularity theory for viscosity solutions, 
e.g., \cite{CIL92}, \cite{CC95},  because of their degenerate dependence  on $\rho$ rather than $\nabla_x\rho$. We refer to the works of Guan \cite{Gu97}, \cite{Gu02} who shows that in certain cases, one is able to "lift" a $C^{1,1}$-regularity of $\rho$ into a statement of $C^\infty$-regularity. Using the existence of smooth viscosity solutions in the uniformly elliptic case where $\bbb(v) \geq \lambda >0$, (\ref{eq:ellent})  could be then justified by the ``vanishing viscosity limit", forming a family of regularized solutions, $\rho^\lambda$ associated with $\bbb^\lambda(s):=\bbb(s)+\lambda I_{d\times d}$ and letting $\lambda\downarrow 
0_+$.  Next comes  the kinetic formulation of (\ref{eq:ell}) which takes the form

\begin{equation}\label{eq:ellkin}
-\nabla_x^\top\cdot\bbb(v)\nabla_x \chi_{\rho}(v) +S(v)\frac{\partial}{\partial v}\chi_{\rho}(v) =\frac{\partial}{\partial v} m(x,v) \quad \hbox{in} \ \dd'(\Gamma\times \Rv),
\end{equation}
for some nonnegative $m\in {\cal  M}^+$ which measures "entropy production". Indeed, for an arbitrary convex ``entropy",
$\eta$, the moments of  (\ref{eq:ellkin}) yield

\begin{eqnarray*}
0 & \geq & -\int \eta''(v)m(x,v)dv = \int \eta'(v)\frac{\partial}{\partial v}m(x,v) dv  = \nonumber \\
&  = & -\sum_{j,k=1}^d \int 
\eta'(v)b_{jk}(v) \frac{\partial^2}{\partial x_j\partial x_k}\chi_\rho(v) dv  +  \int \eta'(v) S(v) 
\frac{\partial}{\partial v}\chi_\rho(v) dv = \nonumber \\
&  = & -\sum_{j,k=1}^d \frac{\partial^2}{\partial x_j\partial x_k} \Betajk(\rho) - \int \big(\eta'(v)S(v)\big)' \chi_\rho(v)dv   
  =  -\trace\Big(\nabla_x\otimes\nabla_x \ \Beta(\rho)\Big) -  \eta'(\rho)S(\rho). \nonumber
\end{eqnarray*}
Thus, the kinetic formulation (\ref{eq:ellkin}) is the dual statement for the entropy inequalities (\ref{eq:ellent}).
We postulate that $\rho$ is a \emph{kinetic solution} of (\ref{eq:ell}) if the corresponding distribution function $\chi_{\rho(x)}(v)$ satisfies (\ref{eq:ellkin}), and we address the regularizing effect of  such kinetic solutions. 

To use the averaging lemma, we first extend (\ref{eq:ellkin})
over the full $\Rdx\times\Rv$-space. Let $\psi$ be  $C_0^\infty(R^+)$-cut-off function,  $\psi(s)\equiv 1$ for $s\geq \epsilon$, and let $\zeta(x)$ denote the smoothed distance function to the boundary, $\zeta(x)=\psi(dist(x, \partial \Gamma))$, then  $f(\vx):=\chi_{\rho(x)}(v)\zeta(x)$ satisfies, in $\dd'(\Rdx\times\Rv)$,
\begin{eqnarray}
\lefteqn{-\sum_{j,k=1}^d b_{jk}(v) \frac{\partial^2}{\partial x_j \partial x_k} f(x,v) =}  \nonumber \\
  & &  \quad \frac{\partial}{\partial v} \zeta(x)m(\vx) + \sum_{j,k=1}^d  \frac{\partial}{\partial x_j} \big(b_{jk}(v) \zeta_{x_k}(x)\chi_\rho(v)\big) + \frac{\partial}{\partial x_k}\big(b_{jk}(v)\zeta_{x_j}(x)\chi_\rho(v)\big)
 + \nonumber \\
  & & \qquad   S(v)\zeta(x)\frac{\partial}{\partial v} \chi_\rho(v)  - \sum_{j,k=1}^d b_{jk}(v)\zeta_{x_j x_k}(x)\chi_\rho(v)  \nonumber \\
& &  \quad =:   \frac{\partial}{\partial v}g_1(\vx)+ \Lambda^\eta_x g_2(\vx) + g_3(\vx) +g_4(\vx). \label{eq:ellprod} 
\end{eqnarray}
 
Assume that (\ref{eq:ell}) is non-degenerate in the sense that there exists
an $\alpha\in (0,1)$ such that

\begin{equation}\label{eq:degell}
 \Big|\Omega_{\bbb}(\delta)\Big| \lesssim \delta^\alpha 
 \ \hbox{and} \  \sup_{|\xi|=1} \sup_{v\in \Omega_{\bbb}(\delta)}|\langle \bbb'(v)\xi,\xi\rangle| \lesssim \delta^{1-\alpha},
\quad  
\Omega_{\bbb}(\delta):= \Big\{v\in I : \ \langle\bbb(v)\xi,\xi\rangle\leq \delta\Big\}.
\end{equation}

We examine the contribution of each of the four terms on the right of (\ref{eq:ellprod}) to the overall $W^{s,1}$-regularity of $\overline{f}$, appealing to the different averaging lemmata term-by-term. The first term on the right involves the bounded measure $g_1=\zeta m$;  averaging lemma \ref{avg} with the usual $(p,q)=(2,1)$ then yields that the corresponding average $\overline{f}_1$ has a $W^{s_1,1}_{loc}$-regularity of order $s_1<2\theta_1, \ \theta_1=\alpha/(3\alpha+2)$. The second term on the RHS of (\ref{eq:ellprod}) involves the gradient of the uniformly bounded term $g_2=\sum b_{jk}(v)\zeta_{x_k}(x)\chi_\rho(v)$; here we can use averaging lemma \ref{thm:homogeneous-averaging} with $\eta=q=1$ to conclude that the corresponding average $\overline{f}_2$ has $W^{s_2,1}_{loc}$-regularity of order $s_2<\theta_2=\alpha/(\alpha+2)$ (in fact, with $q=2$ one concludes a better $W^{s_2,2}$ regularity of order $s_2<\alpha/2$). The remaining two terms on the RHS of (\ref{eq:ellprod}) yield smoother averages and therefore they do not affect the overall regularity dictated by the first two. Indeed,  $\partial_v \chi_\rho(v)$ and hence $g_3(\vx)=S(v)\zeta(x)\partial_v\chi_\rho(v)$ is a bounded measure and  averaging lemma \ref{thm:homogeneous-averaging} with $\eta=N=0$ implies that the corresponding average $\overline{f}_3$ belongs to the smaller Sobolev space, $W^{s_3,1}_{loc}$ of order $s_3<2\theta_3, \ \theta_3=\alpha/(\alpha+2)$.  Finally, the last term on the right of (\ref{eq:ellprod}) consists of the bounded sum, $g_4= - \sum b_{jk}(v)\zeta_{x_j x_k}(x)\chi_\rho$; with $\eta=N=0$ and $q=2$,
 the corresponding average $\overline{f}_4$ has a $W^{s_4,2}_{loc}$ regularity of order $s_4<2\theta_4, \ \theta_4=\alpha/2$. 

Next, we iterate the bootstrap argument we mentioned earlier in the context of hyperbolic conservation laws. The first $W^{s_1,1}$-bound together with the $L^\infty$-bound of $f$ imply a $W^{\sigma_1,2}_{loc}$-bound with $\sigma_1=s_1/2$, which in turn yields the improved regularity of $\overline{f}_1\in W^{s,1}_{loc}, \ s<(1-\theta_1)\sigma_1+2\theta_1$.
Thus, for the first term we can iterate the improved regularity, $s_1\mapsto (1-\theta_1)s_1/2+2\theta_1$, converging to  the same fixed point we had in the parabolic case before, $s_1<2\alpha/(2\alpha+1)$. The second  term   requires a more careful treatment: as we iterate the improved regularity of $\chi_\rho\in W^{s,1}_{loc}$, we can express the  term on  the right of (\ref{eq:ellprod}) as $\Lambda^{\eta_s}g_2$ with $\eta_s:=1-s$ and with $g_2$ standing for  the sum of $L^1$-bounded terms,  $g_2=\Lambda^s\sum b_{jk}(v)\zeta_{x_k}\chi_\rho(v)$. Consequently, averaging lemma \ref{thm:homogeneous-averaging} yields the fixed point iterations $s_2\mapsto (1-\theta_2)s_2/2+(2-\eta_{s_2})\theta_2$ with limiting regularity of order $s_2<\alpha$. The remaining two terms are smoother and do not affect the overall regularity: a similar argument for the third term yields the fixed point iterations, $s_3\mapsto (1-\theta_3)s_3/2+2\theta_3$ with a fixed point $s_3< 2\alpha/(\alpha+1)$, while the fourth term remains in the smaller Sobolev space $W^{\alpha,2}$.
We summarize with the following statement.

\begin{cor} Let $\rho \in L^\infty$ be a kinetic solution of the nonlinear elliptic equation (\ref{eq:ell}) and
assume the non-degeneracy condition (\ref{eq:degell}) holds. Then we have the interior regularity estimate for all $D \subset \Gamma$,
\[
\rho(x) \in W^{s,1}_{loc}(D), \quad s<\left\{ \begin{array}{ll} 
 \alpha,  & {\rm if} \ \ \alpha<1/2, \\ \\
 {\displaystyle \frac{2\alpha}{2\alpha+1}},  & {\rm if} \ \ 1/2 <\alpha <1.
                                     \end{array}\right.
\]
\end{cor}

{\footnotesize

\noindent
\end{document}